\documentclass[preprint,12pt,authoryear]{elsarticle}

\usepackage{lineno}

\journal{}

\usepackage[left=1.5cm, right=1.5cm, top=2cm, bottom=2cm]{geometry} 
\usepackage{multicol}
\usepackage{amsmath, amsfonts, amssymb, amsthm, mathtools, eucal} 
\usepackage{upgreek} 
\usepackage{subfig, tikz, color} 

\theoremstyle{definition}
\newtheorem{remark}{Remark}[section]

\usepackage[bbgreekl]{mathbbol}
\DeclareSymbolFontAlphabet{\mathbb}{AMSb} 
\DeclareMathAlphabet{\bbol}{U}{bbold}{m}{n} 

\usepackage[norelsize,ruled,vlined,titlenumbered,linesnumbered]{algorithm2e}
\SetKwFor{For}{for}{}{end for}
\SetAlgoCaptionSeparator{}
\makeatletter
\newcommand{\nosemic}{\renewcommand{\@endalgocfline}{\relax}}
\newcommand{\dosemic}{\renewcommand{\@endalgocfline}{\algocf@endline}}
\let\oldnl\nl
\newcommand{\nonl}{\renewcommand{\nl}{\let\nl\oldnl}}
\makeatother

\let \SS \relax
\newcommand{\SS}{\mathbb{S}}
\newcommand{\m}[1]{\mathbf{#1}}
\newcommand{\defeq}{\coloneq}
\DeclareMathOperator{\supp}{supp}

\newcommand{\cX}{\mathcal{X}}
\newcommand{\mtau}{\pmb{\uptau}}
\newcommand{\mt}{\m{t}}
\newcommand{\RR}{\mathbb{R}}
\newcommand{\cT}{\mathcal{T}}
\newcommand{\cM}{\mathcal{M}}
\newcommand{\cC}{\mathcal{C}}

\begin{document}
\begin{frontmatter}

\title{An algorithmic approach to direct spline products:\\ procedures and computational aspects} 

\author[TV]{Francesco Patrizi\corref{cor1}} 
\author[UniFi]{Alessandra Sestini}

\affiliation[TV]{organization = {Department of Mathematics, University of Rome Tor Vergata},
            addressline={Via della Ricerca Scientifica 1}, 
            city={Rome},
            postcode={00133}, 
            country={Italy}}

\affiliation[UniFi]{organization = {Department of Mathematics and Computer Science "Ulisse Dini", University of Florence},
            addressline = {Viale Giovanni Battista Morgagni 67/a}, 
            city={Florence},
            postcode={50134}, 
            country={Italy}}

\cortext[cor1]{corresponding author email: patrizi@mat.uniroma2.it}

\begin{abstract}
We introduce an efficient algorithmic procedure for implementing the direct formula that represents the product of splines in the B-spline basis. We first demonstrate the relevance of this direct approach through numerical evidence showing that implicit methods, such as collocation, may fail in some instances due to severe ill-conditioning of the associated system matrices, whereas the direct formula remains robust. We then recast the direct formula into an algorithmic framework based on the Oslo Algorithm and subsequently enhance it, through a factorization of the terms to be computed, to dramatically improve computational efficiency. Extensive numerical experiments illustrate the substantial reduction in computational cost achieved by the proposed method. Implementation aspects are also discussed to ensure numerical stability and applicability.
\end{abstract}




\end{frontmatter}



\section{Introduction}
In this work we present a simple yet effective algorithmic procedure for the direct multiplication of splines. Specifically, we develop an efficient implementation of the direct formula introduced by \citet{morken} for computing the B-spline coefficients of the product of two splines. Despite its theoretical relevance, the formula of \citet{morken} has seen limited practical use, largely due to its cumbersome combinatorial structure and the associated computational costs. As a result, research and applications have often favored alternative approaches, either direct, such as those based on blossoming techniques, as in \citet{cohen, lee, ueda}, or implicit, relying on interpolation or basis transformations, as in \citet{piegl}.

Our contribution revives and significantly enhances the direct approach of \citet{morken}. We recast the original formula into a clear and accessible algorithmic procedure whose core relies on the Oslo Algorithm, developed in \citet{oslo1, oslo2, oslo3}. Building on this formulation, we introduce an improved procedure that dramatically reduces the computational complexity of the direct method. The computation of each B-spline coefficient is reduced to a sequence of matrix–vector products involving rectangular bidiagonal matrices of decreasing size. This structure makes the method transparent, easy to implement, and readily integrable into existing spline-based codes, in contrast to more involved blossoming-based techniques.
Compared to indirect approaches, such as interpolation or change-of-basis methods, the proposed direct procedure is also considerably more numerically stable. In particular, for higher spline degrees, indirect methods often suffer from severe ill-conditioning, leading to substantial loss of accuracy or even numerical failure. Our numerical experiments demonstrate that the direct approach, when implemented with the proposed algorithm, remains robust and yields significantly more accurate results across a wide range of settings.

Beyond its intrinsic theoretical interest, spline multiplication plays a central role in many applications. It is a fundamental building block in the construction of Gram (mass) matrices for $L^2$ projections onto spline spaces and in the assembly of system matrices for Galerkin discretizations, particularly within the framework of isogeometric analysis (IgA, \citet{prefinement1}) and spline-based boundary element methods (IgA-BEM). In IgA, while many applications employ moderate spline degrees, higher-degree regimes are also relevant due to $p$- and $k$-refinement strategies, where increasing the degree provides faster convergence than mesh refinement (see, e.g., the estimates in \cite{prefinement5, prefinement6}) and is routinely used to achieve high accuracy with few degrees of freedom, as in \cite{prefinement2, prefinement3, prefinement4}.
Furthermore, under $p$-refinement, the collocation approach already loses $2$--$3$ digits of accuracy for moderate degrees ($p \leq 10$), independently of mesh fineness. This loss can be significant in high-accuracy simulations and may be further amplified by subsequent ill-conditioned steps, such as the solution of the resulting Galerkin systems.
In this context, the proposed direct approach enhances quadrature schemes and assembly pipelines. Moreover, once the product spline is represented directly in the B-spline basis, integration reduces to the exact integration of B-splines, for which closed-form formulas are available, see, e.g., \citet[Theorem 5]{manni1}. This completely avoids element-wise integration and the use of specialized Gaussian quadrature rules, yielding global (support-wise), exact, and robust integration procedures within spline-based discretization frameworks.
Spline products are also essential in operations involving rational B-splines (NURBS), such as addition, as well as in various geometric and graphical computations, including distance evaluation between spline curves. Moreover, classical spline procedures such as degree elevation can be interpreted as special cases of spline products.

The presentation of the paper is limited to the univariate setting only for clarity. However, the proposed approach extends straightforwardly to multivariate tensor-product spline spaces, thanks to the locality of the algorithm and the separable structure of tensor-product bases. Extensions to multivariate NURBS are likewise immediate, since spline products constitute the core operation in rational arithmetic.

The remainder of the paper is organized as follows. In Section \ref{sec:motivation}, we illustrate the importance of direct spline multiplication in practical applications by comparing the direct approach of \citet{morken}, implemented with our algorithm, to implicit interpolation-based methods. The results reveal accuracy gaps of several orders of magnitude between direct and implicit approaches, with implicit approaches often becoming numerically unreliable. In Section \ref{sec:morken++}, we recall the Oslo Algorithm and the direct formula of \citet{morken}, and reformulate the latter into an algorithmic framework. We analyze the computational challenges of a naive implementation and introduce an improved procedure based on a suitable factorization of the terms, together with additional implementation details to ensure numerical stability. Section \ref{sec:numerics} is devoted to numerical experiments, where we compare the naive and improved implementations of the direct formula and demonstrate the dramatic gains in computational efficiency. Finally, Section \ref{sec:conclusion} concludes the paper and outlines possible directions for future work.

\section{Motivation: a numerical study of stability and correctness}\label{sec:motivation}
In this section we motivate, through numerical experiments, the use of a direct formula for spline products. We show that the standard collocation (interpolation) approach can be unreliable in relevant cases. The issue is the large condition number of the associated matrix, which leads to instability and loss of accuracy. We also show that the improved direct formula of \citet{morken}, described in the next section, provides significantly more accurate results.

We assume familiarity with B-splines and their basic properties. For a brief introduction, see \citet{manni1,manni2}. More comprehensive treatments can be found in \citet{deboor,schumaker}.

Let $p_1, p_2$ be two degrees. Let $\pmb{\uptau}^1 = (\tau_j^1)_{j=1}^{n_1 + p_1 + 1}$ and $\pmb{\uptau}^2 = (\tau_j^2)_{j=1}^{n_2 + p_2 + 1}$ be two knot vectors, assumed open with respect to $p_1$ and $p_2$ (that is, with first and last knots of multiplicity $p_j + 1$, for $j = 1, 2$, respectively). Denote by $\SS_{\pmb{\uptau}^1}^{p_1}$ and $\SS_{\pmb{\uptau}^2}^{p_2}$ the corresponding spline spaces. Then $\dim \SS_{\pmb{\uptau}^j}^{p_j} = n^j$ for $j=1,2$.

Let $f \in \SS_{\pmb{\uptau}^1}^{p_1}$ and $g \in \SS_{\pmb{\uptau}^2}^{p_2}$. Their product $h \defeq f \cdot g$ belongs to $\SS_{\m{t}}^p$, where $p = p_1 + p_2$ and $\m{t} = (t_j)_{j=1}^{m + p + 1}$ is open with respect to $p$. The breakpoints of $\m{t}$ are the distinct values in $\pmb{\uptau}^1 \cup \pmb{\uptau}^2$. The knot multiplicities are chosen to preserve the lowest smoothness of the factors:
\begin{equation}\label{eq:tmultiplicities}
\mu(t_i) = \left\{\begin{array}{ll}
\max\{ p_1 + \mu^2(t_i), \, p_2 + \mu^1(t_i)\} & \text{if $\mu^1(t_i) > 0$ and $\mu^2(t_i) > 0$},\\\\
p_1 + \mu^2(t_i) & \text{if $\mu^1(t_i) = 0$ and $\mu^2(t_i) > 0$},\\\\
p_2 + \mu^1(t_i) & \text{if $\mu^1(t_i) > 0$ and $\mu^2(t_i) = 0$},
\end{array}\right.
\end{equation}
where $\mu^j(t_i)$ is the multiplicity of $t_i$ in $\pmb{\uptau}^j$.

The goal is to determine the coefficients $\m{b} \defeq (b_i)_{i=1}^m$ such that
\begin{equation}\label{eq:hBspline}
h(x) = \sum_{i=1}^m b_i B_{i,\m{t}}(x).
\end{equation}

A widespread approach is collocation. This is an implicit method, since it requires solving a linear system. It is often used because of its simplicity. The idea is to interpolate the values of $f \cdot g$.
Since $\dim \SS_{\m{t}}^p = m$, we impose
$$
h(x_i) = f(x_i) g(x_i) \eqqcolon h_i, \quad i=1,\dots,m,
$$
at $m$ distinct nodes $\{x_i\}_{i=1}^m$. This yields the linear system $\bbol{C}\m{b} = \m{h}$, where $\bbol{C}$ is the B-spline collocation matrix
$$
\bbol{C} \defeq \begin{bmatrix}
B_{1,\m{t}}(x_1) & \cdots & B_{m,\m{t}}(x_1)\\
\vdots & \ddots & \vdots\\
B_{1,\m{t}}(x_m) & \cdots & B_{m,\m{t}}(x_m)
\end{bmatrix}.
$$
Such a matrix is nonsingular if and only if the nesting conditions of \cite{schoenbergwhitney} hold, i.e., $B_{i,\m{t}}(x_i) > 0$ for all $i$. In that case $\bbol{C}$ is in fact totally positive. A common choice to verify the conditions is to take the Greville abscissae as nodes:
$$
x_i \defeq \frac{t_{i+1} + \cdots + t_{i+p}}{p}, \quad i=1,\dots,m.
$$
With a standard ordering, $\bbol{C}$ is banded with bandwidth $p+1$. The cost of LU factorization is therefore $\frac{2}{3} m p^2$.

Alternative implicit approaches use B\'ezier extraction, see \citet{piegl}. These involve knot insertion, Bernstein multiplication, and knot removal. However, conversions to Bernstein form are ill-conditioned processes, see \citet[Section 6.5]{rida}. Direct approaches based on blossoming have also been proposed, see \cite{ueda,lee}. These are equivalent to the formula of \citet{morken}, but are more involved.

However, in the remainder of this section we show, through numerical tests, that direct formulas are sometimes necessary. We compare the accuracy achieved with our enhanced version of the direct approach of \citet{morken} and collocation. In the tests, splines are defined on open knot vectors with simple internal knots, unless stated otherwise (as in, e.g., Section \ref{sec:galerkin}). We use uniform knots. Experiments with non-uniform knots, including highly clustered ones, showed no qualitative differences and are omitted.
Errors are measured using a discrete relative $L^\infty$ norm, defined as follows. Let $\cX \defeq \{x_j\}_{j=1}^N$ be a set of uniform points. For the direct method (denoted by $M$) and collocation (denoted by $C$), we define
$$
e^M(x) \defeq \frac{|h^M(x) - h(x)|}{\max_{\cX} |h(x)|}, \quad
e^C(x) \defeq \frac{|h^C(x) - h(x)|}{\max_{\cX} |h(x)|}.
$$
In all tests we fix $N = 201$. For each test, we present the errors $\max_{\cX} e^M$ and $\max_{\cX} e^C$ under varying degrees or mesh sizes, and the estimated condition number of the matrix $\bbol{C}$, using the command \texttt{condest} of Matlab. Across all cases, $e^M$ remains at or near machine precision (between \texttt{1e-16} and \texttt{1e-14}), even for high polynomial degrees. In contrast, $e^C$ deteriorates, often rapidly, with increasing degree, in line with the growth of its condition number, and may lead to numerical breakdown.

\subsection{Product of a fixed B-spline, or spline, and polynomials of increasing degrees}\label{sec:splinepoly}
In this test we consider the product of a fixed spline with polynomials of increasing degree, a setting arising, e.g., in B-spline weighted quadrature in IgA-BEM, where spline basis functions are multiplied by polynomial approximations of boundary kernels (see \citet{stokes}). In such contexts, $p$-refinement on the polynomial factor improves accuracy, and spline--polynomial products are required prior to integration when assembling the discrete system.

\begin{figure}
\centering
\subfloat[]{\includegraphics[width = 0.4\textwidth, page = 1]{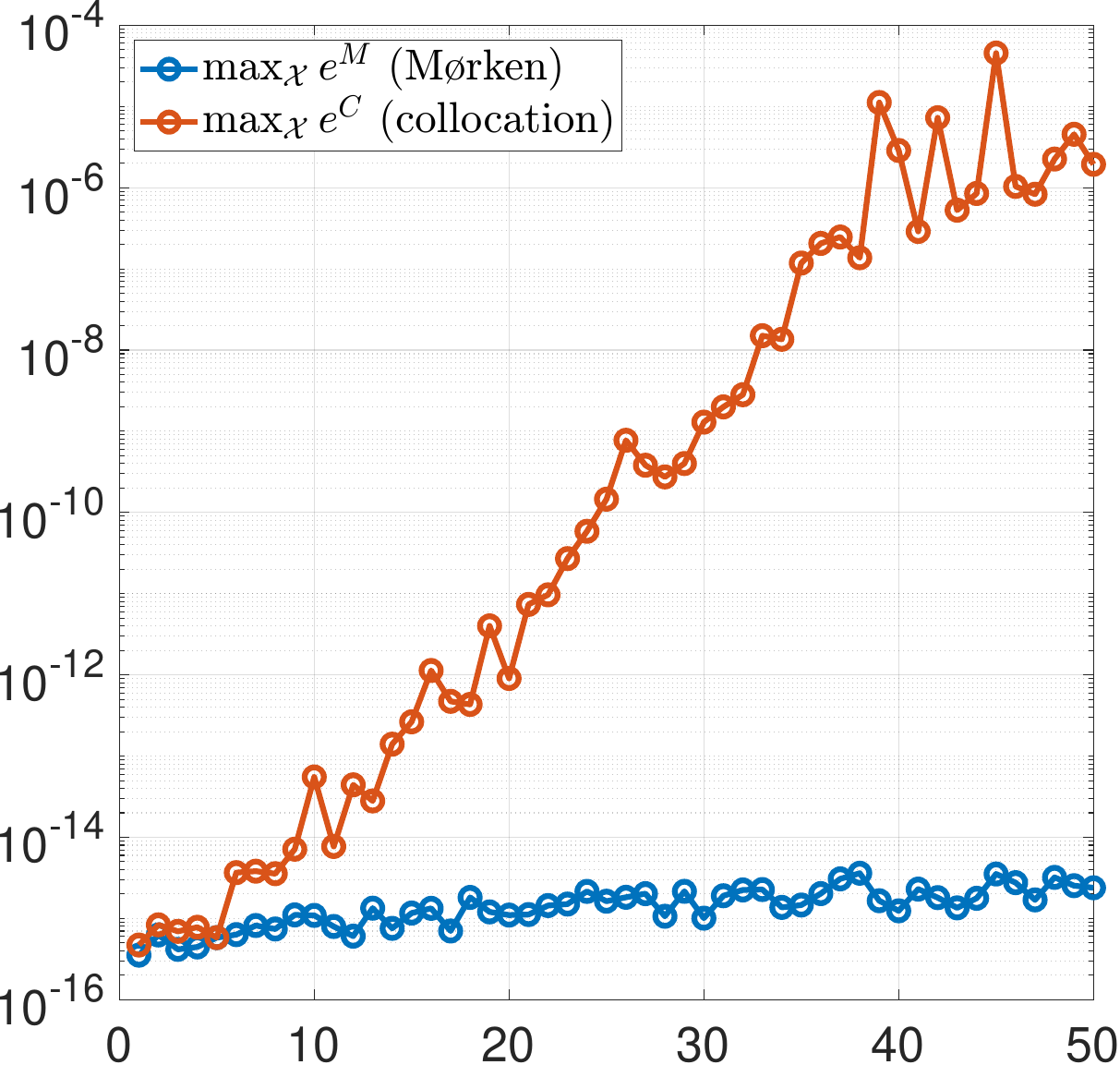}}\qquad
\subfloat[]{\includegraphics[width = 0.4\textwidth, page = 2]{tests}}
\caption{Error comparison and condition number of $\bbol{C}$ for the product of a fixed cubic B-spline and a polynomial of increasing degree. (a) $\max_{\cX} e^M$ and $\max_{\cX} e^C$ vs.\ polynomial degree. (b) Condition number of $\bbol{C}$.}
\label{fig:test1}
\end{figure}

We fix $f$ as a B-spline of degree $p_1 = 3$ on $[0,1]$ and let $g$ be a polynomial on $[0,1]$ with random coefficients in $[-1,1]$ and degree $p_2 \in \{1,\dots,50\}$. Results are shown in Figure \ref{fig:test1} (a)--(b).
We also consider $f$ as a generic spline of degree $p_1 = 3$ on $[0, 1]$ and random coefficients in $[-1,1]$, with $g$ as above. The results are consistent with the single B-spline case.

\subsection{The Galerkin isogeometric case with $p$- and $k$-refinements}\label{sec:galerkin}
In this test we consider the product of two B-splines defined on the same uniform knot vector $\pmb{\uptau}$ on $[0,1]$, with equal degree $p_1 = p_2 \in \{3,\dots,50\}$. As the degree increases, knot multiplicities are adjusted to preserve global $C^2$ smoothness, as in standard $p$-refinement procedures for Galerkin discretizations (see \citet{prefinement1,prefinement2}).

\begin{figure}
\centering
\subfloat[]{\includegraphics[width = 0.4\textwidth, page = 3]{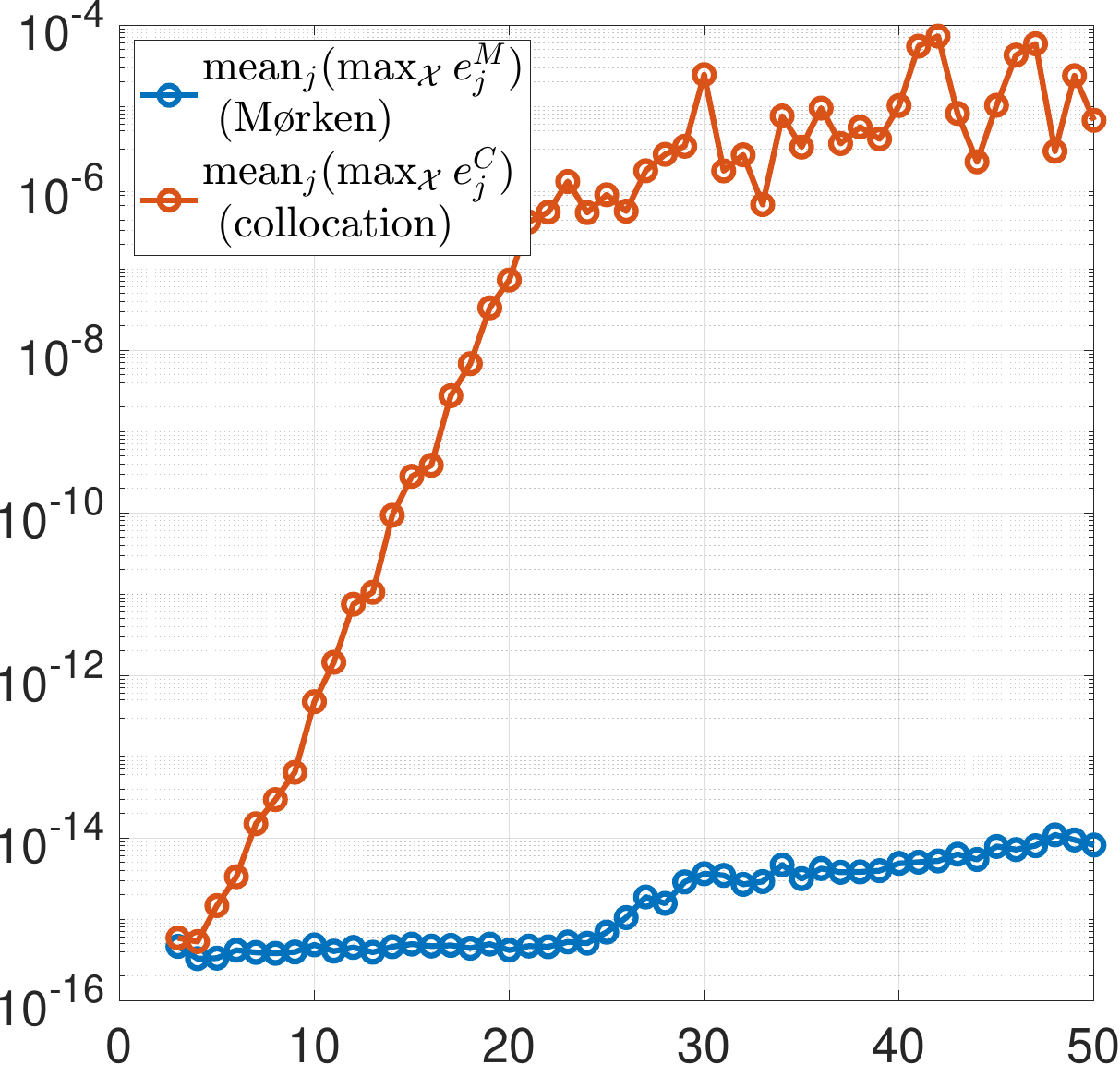}}\quad
\subfloat[]{\includegraphics[width = 0.4\textwidth, page = 4]{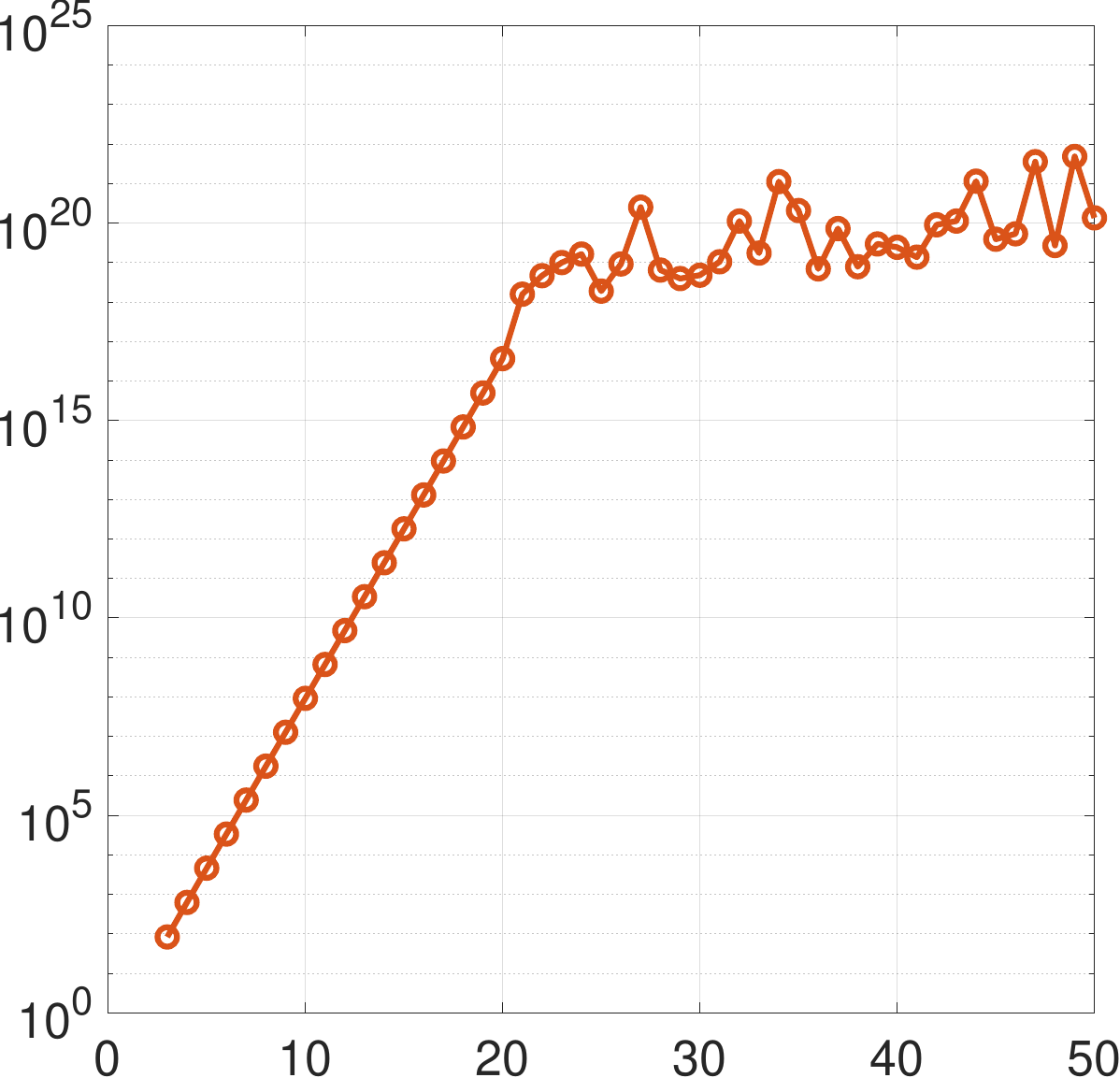}}
\caption{Error comparison and collocation condition number for $B_{i,\pmb{\uptau}} \cdot B_{j,\pmb{\uptau}}$ on a common knot vector $\pmb{\uptau}$ with fixed smoothness ($C^2$) and increasing degree. For fixed $i$, mean of $\max_\cX e_j^M$ and $\max_\cX e_j^C$ vs.\ polynomial degree (a) and mean condition number of $\bbol{C}_j$ (b), over all $j$ such that $B_i$ and $B_j$ have overlapping supports.}
\label{fig:test2}
\end{figure}

\begin{figure}
\centering
\subfloat[]{\includegraphics[width = 0.4\textwidth, page = 5]{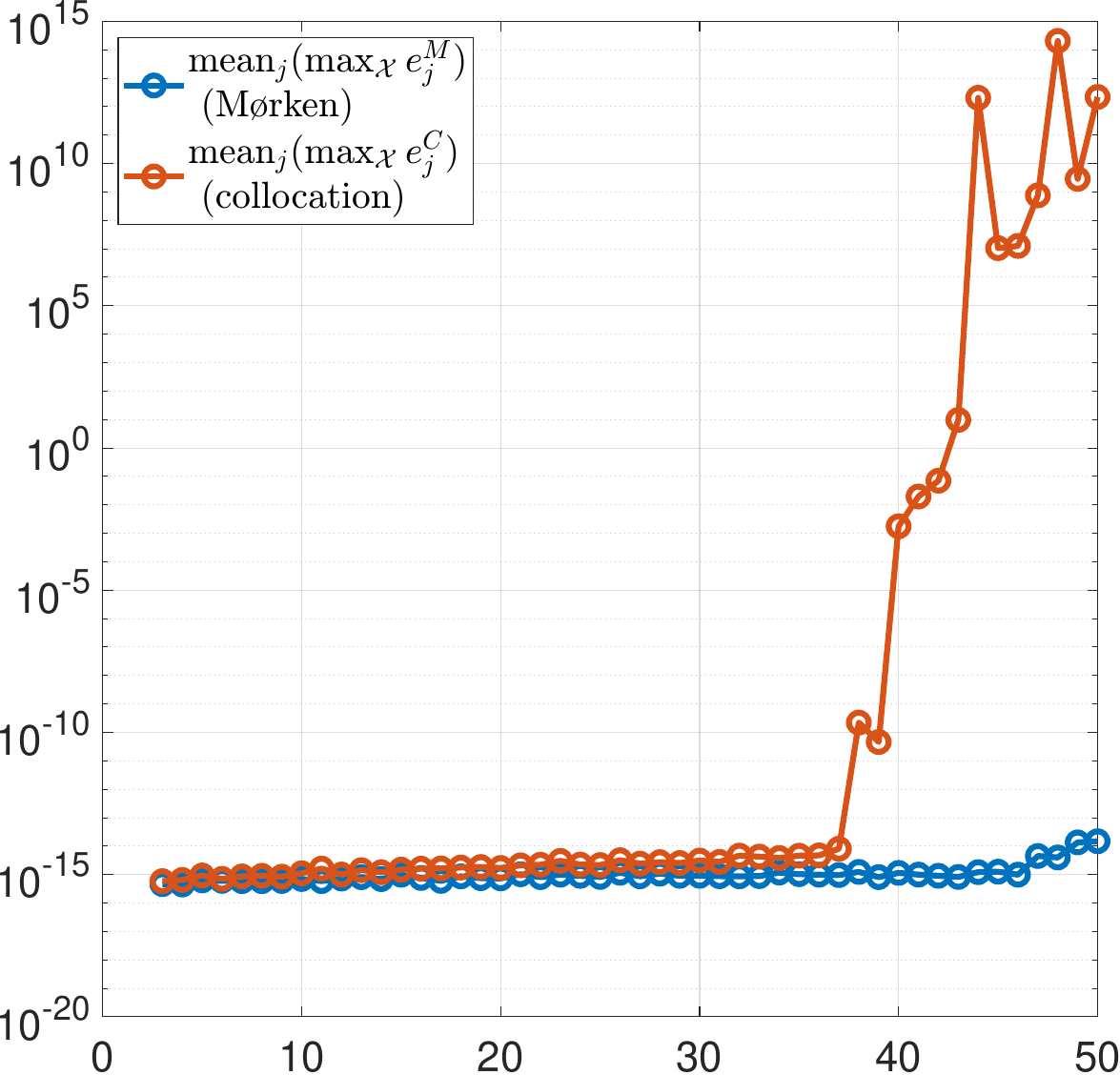}}\quad
\subfloat[]{\includegraphics[width = 0.4\textwidth, page = 6]{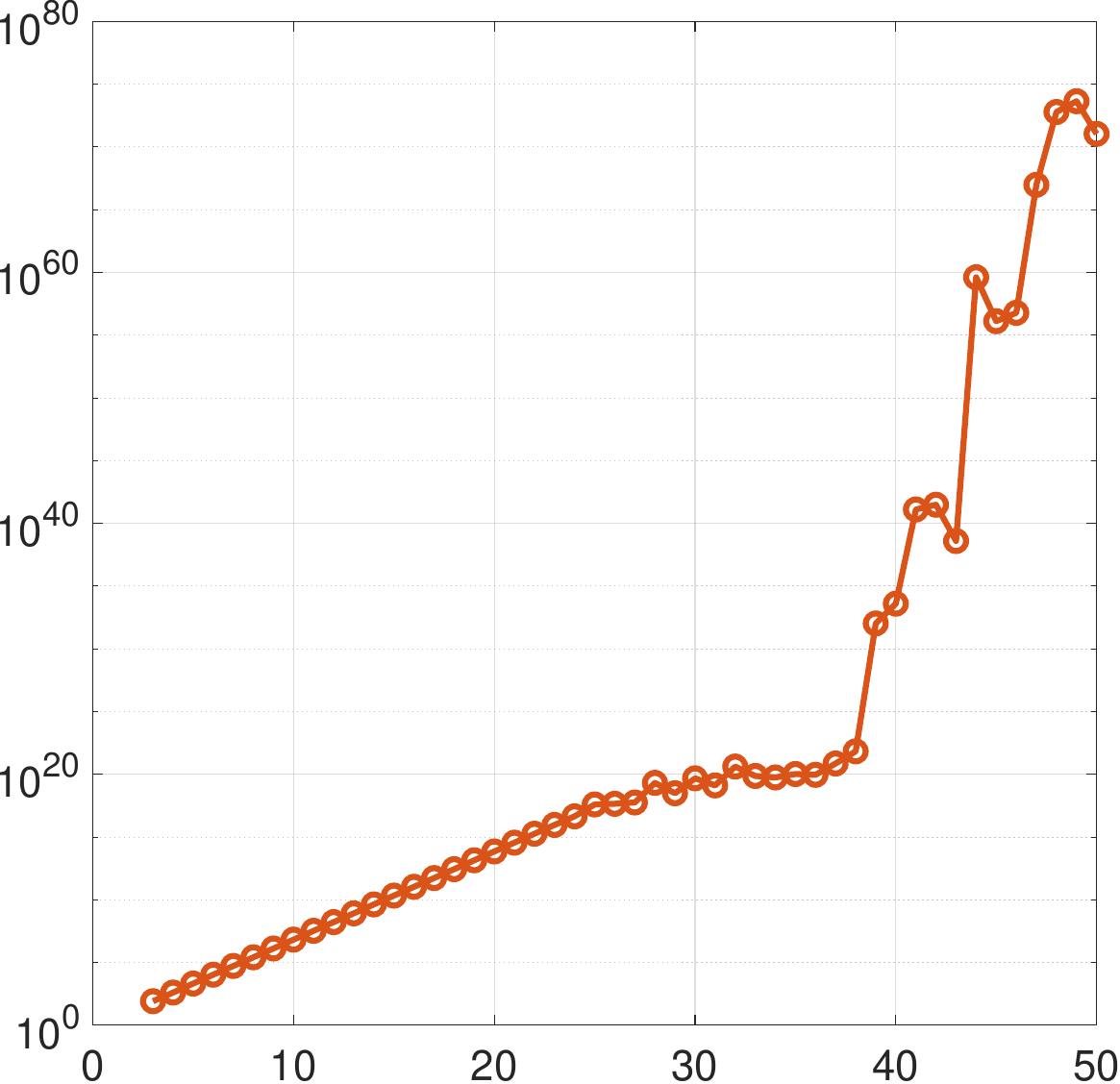}}
\caption{Error comparison and collocation condition number for $B_{i,\pmb{\uptau}} \cdot B_{j,\pmb{\uptau}}$ with maximal smoothness and increasing degree. For fixed $i$, mean of $\max_\cX e_j^M$ and $\max_\cX e_j^C$ vs.\ polynomial degree (a) and mean condition number of $\bbol{C}_j$ (b), over all $j$ such that $B_i$ and $B_j$ have overlapping supports.}
\label{fig:test21}
\end{figure}

For each degree, we fix $f = B_{i,\pmb{\uptau}}$ and consider all B-splines $B_{j,\pmb{\uptau}}$ with overlapping support, i.e., $\supp B_{i,\pmb{\uptau}} \cap \supp B_{j,\pmb{\uptau}} \neq \varnothing.$
The reported error is the average of the relative maximal errors over all such pairs. Results are shown in Figure \ref{fig:test2} (a)--(b).
We also examine the $k$-refinement regime, where smoothness increases to $C^{p_j-1}$ with the degree. This represents the most challenging case for collocation. Results are reported in Figure \ref{fig:test21} (a)--(b).

\subsection{Product of splines on same knot vector and same degree}
\begin{figure}
\centering
\subfloat[]{\includegraphics[width = 0.4\textwidth, page = 7]{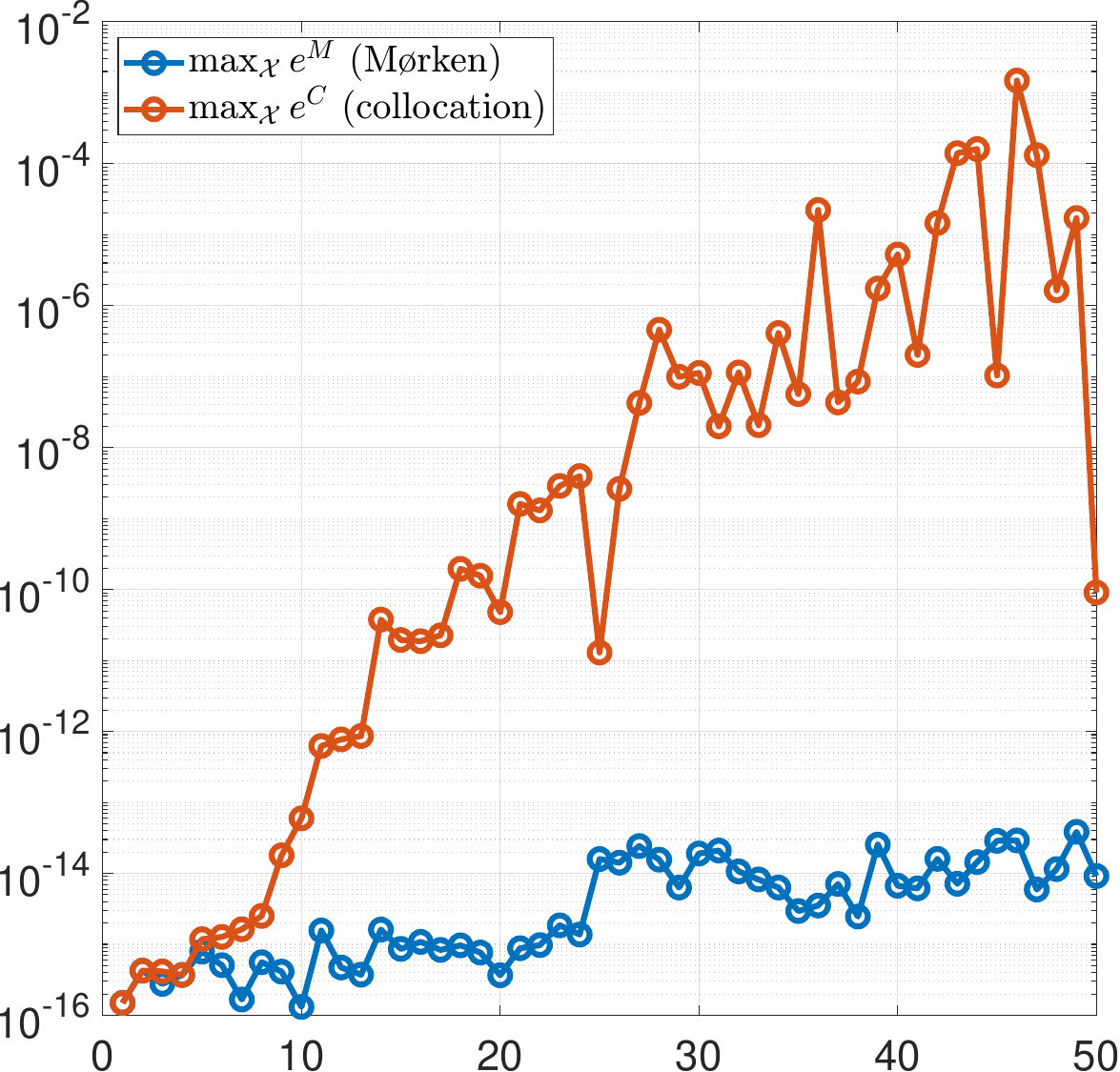}}\qquad
\subfloat[]{\includegraphics[width = 0.4\textwidth, page = 8]{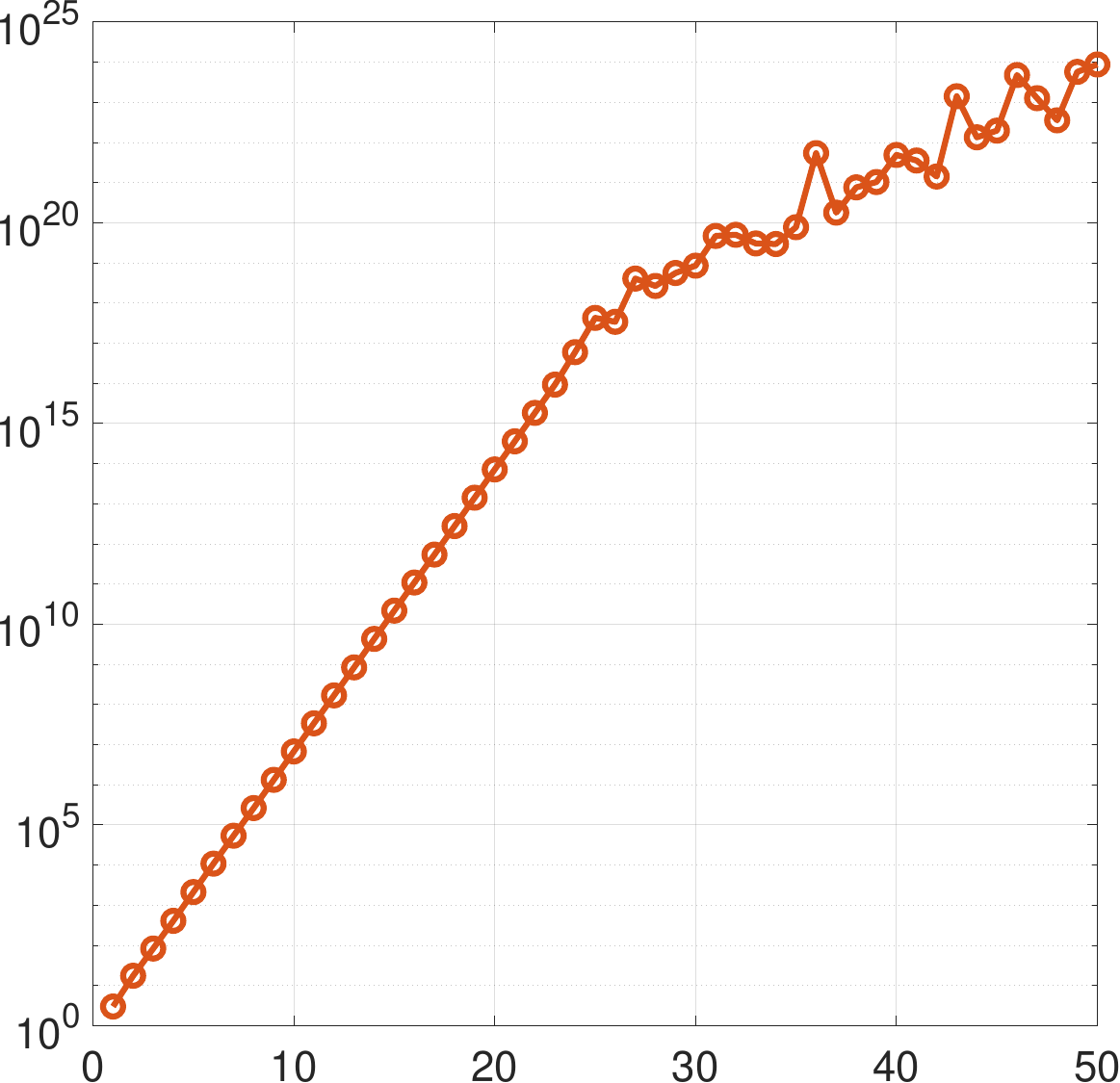}}
\caption{Error comparison and condition number of $\bbol{C}$ for the product of splines on a common knot vector $\pmb{\uptau}$ with increasing degree. (a) $\max_\cX e^M$ and $\max_\cX e^C$ vs. polynomial degree. (b) Estimated condition number of $\bbol{C}$.}
\label{fig:test3}
\end{figure}

As a further test, we consider the product of two splines defined on the same knot vector $\pmb{\uptau}$ on $[0,1]$, with equal degree $p_1 = p_2 \in \{1,\dots,50\}$. This generalizes the previous setting (Galerkin case with $k$-refinement) by replacing B-splines with generic splines.
We take $\pmb{\uptau}$ with $5$ uniform breakpoints, and define $f$ and $g$ as splines with random coefficients in $[-1,1]$. Results are reported in Figure \ref{fig:test3} (a)--(b).

\subsection{Product of a fixed spline with splines on increasingly finer meshes}
In this test we consider the product of a fixed spline with splines defined on successive refinements of its knot vector. This setting arises, e.g., in IgA-BEM quadrature schemes \cite{helm}, where a spline is multiplied by a spline quasi-interpolant defined on a refined mesh.

We fix $f$ as a spline of degree $p_1 = 3$ on a knot vector $\mtau^1$ with $5$ uniform breakpoints in $[0,1]$, with random coefficients in $[-1,1]$. The second factor $g$ is a spline of fixed degree $p_2$ with random coefficients in $[-1,1]$, defined on knot vectors $\mtau^2$ obtained by dyadic refinement of $\mtau^1$, with $2^n + 3$ breakpoints, $n = 1,\dots,10$.

\begin{figure}
\centering
\subfloat[]{\includegraphics[width = 0.45\textwidth, page = 9]{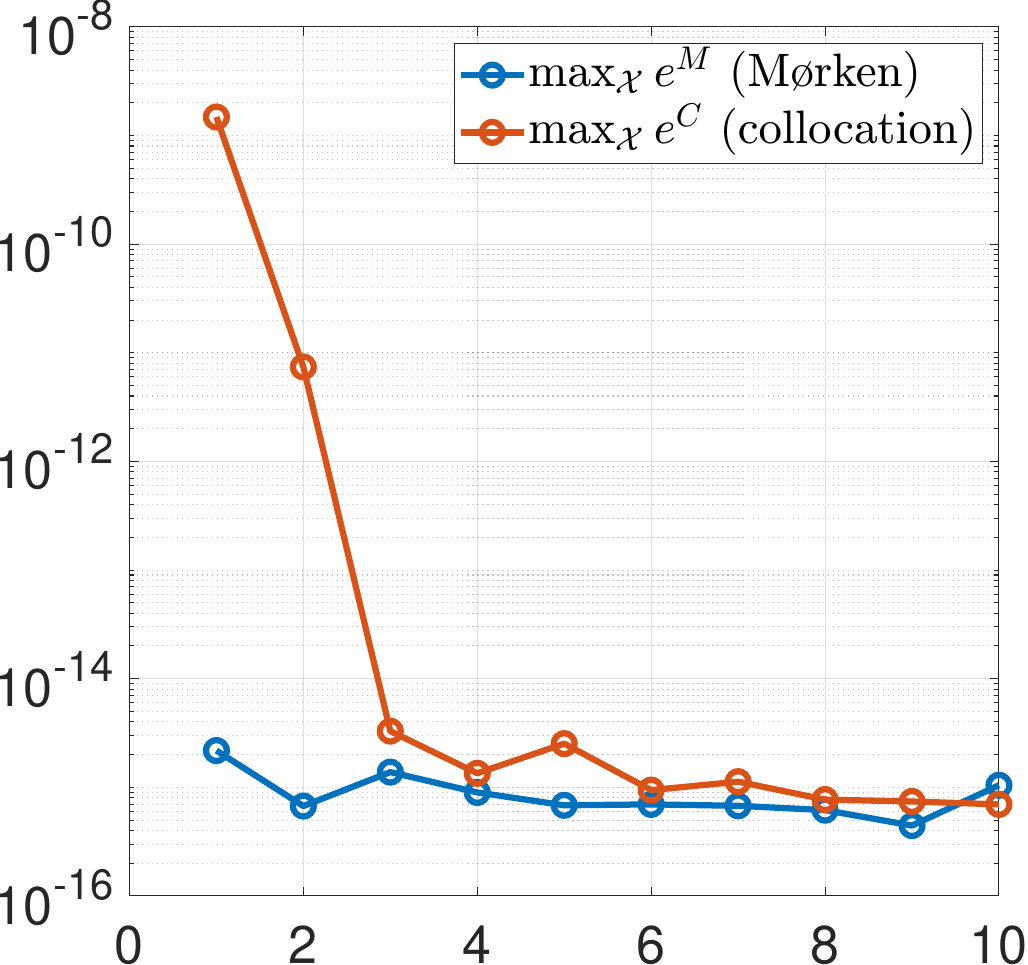}}\quad
\subfloat[]{\includegraphics[width = 0.45\textwidth, page = 10]{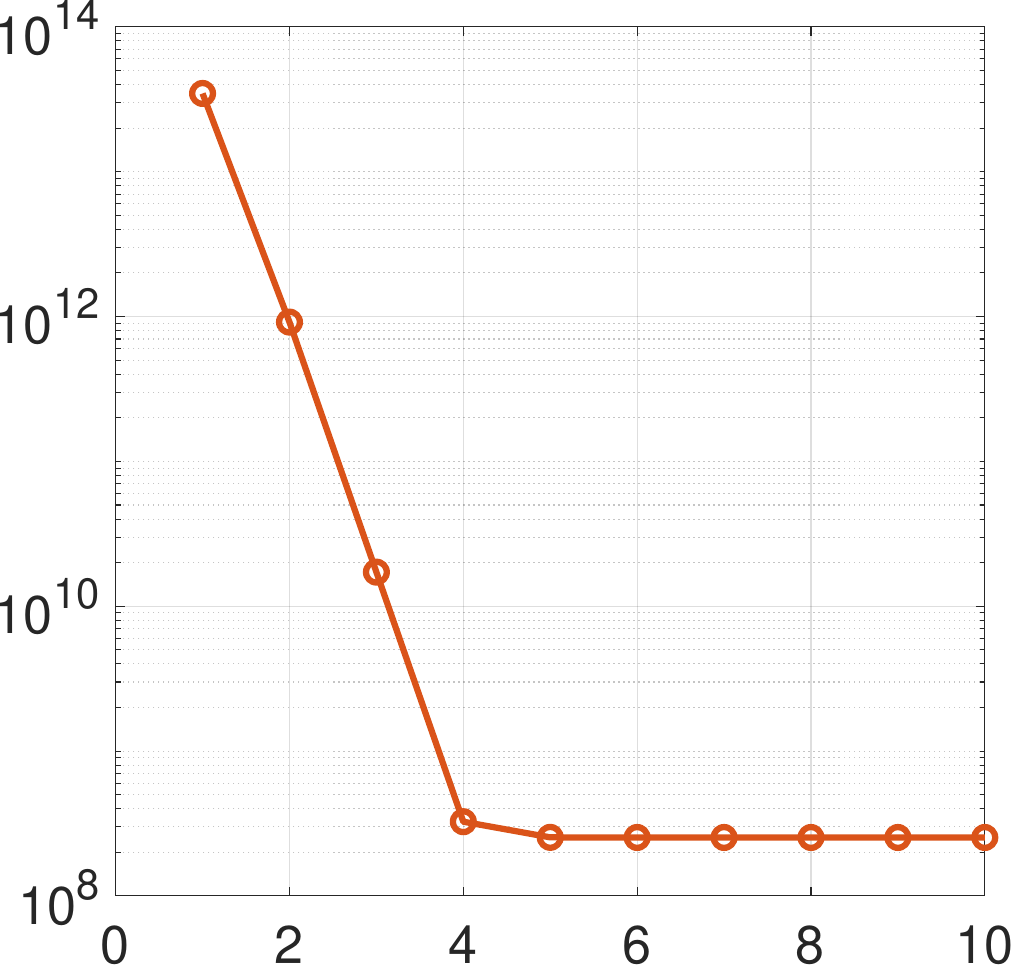}}
\caption{Error comparison and condition number of $\bbol{C}$ for the product of a cubic spline and a spline of degree $30$. The first spline is defined on a fixed knot vector $\mtau^1$ with $5$ uniform breakpoints, the second on $\mtau^2$ with $2^n+3$ breakpoints, $n=1,\ldots,10$. (a) $\max_\cX e^M$ and $\max_\cX e^C$ vs.\ $n$. (b) Estimated condition number of $\bbol{C}$.}
\label{fig:test41}
\end{figure}

For low degrees (e.g., $p_2 = 3$), both approaches yield errors of machine precision across all meshes. For higher degrees, the behavior changes: Figure \ref{fig:test41} (a)--(b) shows the case $p_2 = 30$, where collocation exhibits non-negligible errors on coarse meshes before aligning with the direct approach under refinement. However, the refinement level at which accuracy is recovered is not known \emph{a priori}, and can only be inferred indirectly, e.g., from the condition number of $\bbol{C}$. 

\section{Improving the direct formula}\label{sec:morken++}
In this section we review the direct formula of \citet{morken} and propose a refined variant. This new formulation reduces the computational effort and makes the approach practical. 

The original presentation is given as a single formula. We instead express it as a simple procedure based on the Oslo Algorithm. This algorithm is central both to the original construction and to our improvement. We therefore begin with a brief review of the Oslo Algorithm. We then describe the method of \citet{morken} from an algorithmic viewpoint and introduce our enhanced procedure.

\subsection{The Oslo Algorithm}

Let $p \geq 1$. Let $\pmb{\uptau} = (\tau_j)_{j=1}^{n+p+1}$ and $\m{t} = (t_i)_{i=1}^{m+p+1}$ be two open knot vectors of degree $p$, defined on the same interval $I \subseteq \RR$. Assume $\pmb{\uptau} \subseteq \m{t}$, counting multiplicities. Then $\SS_{\pmb{\uptau}}^p \subseteq \SS_{\m{t}}^p$.

Let $q \in \SS_{\pmb{\uptau}}^p$ and let $\m{c} = (c_j)_{j=1}^n$ be its B-spline coefficients. The Oslo Algorithm computes the coefficients $\m{b} = (b_i)_{i=1}^m$ of $q$ in the B-spline basis of $\SS_{\m{t}}^p$. This corresponds to performing multiple knot insertions at once.

We adopt the matrix formulation of \citet[Section 4.2.3]{oslo4}. For $d = 1,\dots,p$, $k \in \{p+1,\dots,n\}$, and $t \in I$, define the bidiagonal $d \times (d+1)$ matrices $\bbol{R}_d^k(t)$ by
$$
\bbol{R}_d^k(t) \defeq \begin{bmatrix}
\frac{\tau_{k + 1} - t}{\tau_{k + 1} - \tau_{k + 1 - d}} & \frac{t - \tau_{k + 1 - d}}{\tau_{k + 1} - \tau_{k + 1 - d}} & 0 & \cdots & 0\\
0 & \frac{\tau_{k + 2} - t}{\tau_{k + 2} - \tau_{k + 2 - d}} & \frac{t - \tau_{k + 2 - d}}{\tau_{k + 2} - \tau_{k + 2 - d}} & \cdots & 0\\
\vdots & \vdots & \ddots & \ddots & \vdots \\
0 & 0 & \cdots & \frac{\tau_{k + d} - t}{\tau_{k + d} - \tau_k} & \frac{t - \tau_k}{\tau_{k + d} - \tau_k}
\end{bmatrix},
$$
with the convention that zero denominators yield zero entries.

The algorithm is summarized in Algorithm \ref{alg:oslo}. For each coefficient $b_i$, select the index $k$ such that $t_i \in [\tau_k, \tau_{k+1})$. This determines the local coefficient vector $\m{c}^k \defeq (c_j)_{j=k-p}^k$ and the local knot vector $\pmb{\uptau}^k \defeq (\tau_j)_{j=k+1-p}^{k+p}$. Let also $\mt^i \defeq (t_j)_{j=i+1}^{i+p}$.
Using these data, construct the matrices $\bbol{R}_d^k(t_{i+d})$ for $d=1,\dots,p$. The algorithm then applies successive matrix--vector products. The matrix sizes decrease from $p \times (p+1)$ to $1 \times 2$. Accordingly, the vector $\m{c}^k$ is reduced step by step, from length $p+1$ to a scalar. This scalar is the coefficient $b_i$.
\begin{algorithm}[h!]
\caption{Oslo Algorithm - $\m{b} = \texttt{oslo\_algorithm}(p, \mtau, \m{c}, \m{t})$}\label{alg:oslo}\vspace{0.1cm}
\nonl\begin{minipage}{\textwidth}
\rule{0.95\textwidth}{0.025cm}

\textbf{Input:}\begin{itemize}\vspace{-0.25cm}\setlength{\itemindent}{-0.2cm}
\item[$p\colon$] degree\vspace{-0.05cm}
\item[$\mtau\colon$] $(\tau_j)_{j=1}^{n + p + 1}$ open knot vector on $I$\vspace{-0.05cm}
\item[$\m{c}\colon$] B-spline coefficients in $\SS_{\mtau}^p$\vspace{-0.05cm}
\item[$\m{t}\colon$] $(t_j)_{j=1}^{m + p + 1}$ open knot vector on $I$ and such that $\mtau \subseteq \mt$\vspace{-0.05cm}
\end{itemize}
\noindent\textbf{Output:}\begin{itemize}\vspace{-0.25cm}\setlength{\itemindent}{-0.2cm}
\item[$\m{b}\colon$] $(b_i)_{i=1}^m$ B-spline coefficients in $\SS_\m{t}^p$\vspace{-0.5cm}
\end{itemize}
\rule{0.95\textwidth}{0.025cm}\vspace{-0.1cm}\end{minipage}

\linespread{1.35}\selectfont
\For{$i = 1, \ldots, m$}{
Find $k$ such that $t_i \in [\tau_k, \tau_{k+1})$\;
Set $\m{c}^k\defeq (c_{k-p}, \ldots, c_k)^T \subseteq \m{c}$\;
Assign $\displaystyle b_i \defeq \bbol{R}_1^k(t_{i + 1})\bbol{R}_2^k(t_{i + 2})\cdots\bbol{R}_p^k(t_{i + p})\m{c}^k$\;}
\end{algorithm}

Thanks to the bidiagonal structure of $\bbol{R}_d^k$, explicit assembly of these matrices is not required. A matrix-free implementation is possible via a de Boor-type procedure, see \citet[Section 2.5]{oslo4}, as schematized in Algorithm \ref{alg:deboor}. In practice, the matrix multiplications in the last step of Algorithm \ref{alg:oslo} are replaced by $b_i = \texttt{deboor}(p, \mtau^k, \m{c}^k, \m{t}^i)$.
\begin{algorithm}[h!]
\caption{de Boor--like Algorithm - $b_i = \texttt{deboor}(p, \mtau^k, \m{c}^k, \m{t}^i)$}\label{alg:deboor}
\nonl\begin{minipage}{\textwidth}
\rule{0.95\textwidth}{0.025cm}

\textbf{Input:}\begin{itemize}\vspace{-0.25cm}\setlength{\itemindent}{-0.2cm}
\item[$p\colon$] degree\vspace{-0.05cm}
\item[$\mtau^k\colon$] $(\tau_j)_{j=k+1-p}^{k + p} \subseteq \mtau$ knot subvector needed to compute $b_i$\vspace{-0.05cm}
\item[$\m{c}^k\colon$] $(c_j^k)_{j=1}^{p+1} \defeq (c_j)_{j=k-p}^k\subseteq \m{c}$ subset of the B-spline coefficients in $\SS_{\mtau}^p$\vspace{-0.05cm}
\item[$\m{t}^i\colon$] $(t_j)_{j=i+1}^{i + p}\subseteq \mt$ knot sub-vector needed to compute $b_i$\vspace{-0.05cm}
\end{itemize}
\noindent\textbf{Output:}\begin{itemize}\vspace{-0.25cm}\setlength{\itemindent}{-0.2cm}
\item[$b_i\colon$] $i$th B-spline coefficient in $\SS_\m{t}^p$\vspace{-0.1cm}
\end{itemize}
\textbf{Remark:} all vector products and divisions are meant element-wise.\hfill\vspace{-0.1cm}
\rule{0.95\textwidth}{0.025cm}\vspace{0.1cm}\end{minipage}

\linespread{1.35}\selectfont
\For{$d = p, p -1, \ldots, 1$}{
Extract $\mtau_d^{k,1} \defeq (\tau_j)_{j=k+1-d}^k$ and $\mtau_d^{k,2} \defeq (\tau_j)_{j=k+1}^{k+d}$ from $\mtau^k$\;
Set $\pmb{\upomega}_d \defeq (t_{i + d} - \mtau_d^{k,1}) / (\mtau_d^{k, 2} - \mtau_d^{k, 1})$\;
Define $\m{c}^{k, 1} \defeq \m{c}^k \setminus\{c_1^k\}$ and $\m{c}^{k, 2} \defeq \m{c}^k \setminus \{c_{d+1}^k\}$\;
Update $\m{c}^k \leftarrow (1 - \pmb{\upomega}_d)\m{c}^{k, 1} + \pmb{\upomega}_d\m{c}^{k, 2}$\;
}
Set $b_i = c^k$\;
\end{algorithm}

We conclude with remarks on implementation and computational aspects.

\begin{remark}\label{oss:notopen}
The Oslo Algorithm requires $\mtau$ and $\mt$ to be open with respect to $p$. This is not restrictive. If $\mtau$ is not open, we extend it to an open knot vector $\widehat{\mtau}$ by repeating the first and last knots. The additional B-spline coefficients are set to zero. The same applies to $\mt$, yielding $\widehat{\mt}$. After computing $\m{b}$ on $\widehat{\mt}$, we retain only the coefficients corresponding to the original knot vector $\mt$.
\end{remark}

\begin{remark}
The cost of multiplying a $d \times (d+1)$ bidiagonal matrix by a vector of length $d+1$ is $3d$ flops, as each row requires two multiplications and one addition. For each coefficient $b_i$, the algorithm performs $p$ such products, with decreasing sizes. The total cost per coefficient is
$$
\sum_{d=1}^p 3d = 3\frac{p(p+1)}{2} \doteq \frac{3}{2}p^2,
$$
where $\doteq$ denotes equality up to lower-order terms. This estimate excludes matrix assembly (which is consistent with the estimate for the collocation approach, where we consider only the LU factorization cost).

The total cost of the Oslo Algorithm is therefore $\frac{3}{2}mp^2$. However, each coefficient $b_i$ is computed independently. The method is thus well suited for parallelization.
\end{remark}

\begin{remark}\label{oss:localvec}
The computation of each $b_i$ is local. It depends only on the subvector $\mt^i$ of $\mt$ and on the local data $\m{c}^k$ and $\mtau^k$. As a consequence, $\mt$ need not be a global refinement of $\mtau$. It is sufficient that $\mt$ locally refines $\mtau$. In particular, $[t_1, t_{m+p+1}]$ may be a subinterval of $[\tau_1, \tau_{n+p+1}]$. In this case, the Oslo Algorithm yields the coefficients of the restriction of $q \in \SS_{\mtau}^p$ to $[t_1, t_{m+p+1}]$ in $\SS_{\mt}^p$.
\end{remark}

\begin{remark}
Algorithm \ref{alg:deboor} generalizes the de Boor algorithm for spline evaluation: if $\mt^i$ is replaced by $(t,\ldots,t)$, the output $b_i$ is the value of the spline at $t$. 
The Oslo Algorithm also generalizes B\'ezier extraction. If $\mt$ has the same breakpoints as $\mtau$, each with multiplicity $p+1$, the algorithm yields the Bernstein representation.
\end{remark}

\begin{remark}
The product $\bbol{R}_1^k(t_{i+1}) \cdots \bbol{R}_p^k(t_{i+p}) \eqqcolon \pmb{\upalpha}^T$ defines a row vector of length $p+1$. Its entries are called discrete B-splines of order $p+1$,  typically denoted as $\alpha_{j,p+1,\mtau,\mt}(i)$. The coefficient $b_i$ can be written as $b_i = \pmb{\upalpha}^T \m{c}^k$.
Precomputing $\pmb{\upalpha}$ may be advantageous when expressing many splines of $\SS_{\mtau}^p$ in $\SS_{\mt}^p$. However, this is not the case for the direct formula of \citet{morken}. We therefore avoid forming $\pmb{\upalpha}$. Instead, we update $\m{c}^k$ after each multiplication. These operations are performed in a matrix-free fashion via a de Boor-type procedure. This avoids both matrix assembly and storage.
The original formulation of \citet{morken} is instead written in terms of these discrete B-splines and presented as a single formula.
\end{remark}

\subsection{The M{\o}rken formula}
The coefficient $b_i$ in \eqref{eq:hBspline} for the product $h = f \cdot g$ is given by \citet[Theorem 3.1]{morken} as
\begin{equation}\label{eq:diequation}
b_i = \frac{1}{{p \choose p_1}}\sum_{P \in \mathcal{C}}\left[(\pmb{\upalpha}^P)^T\m{c}^{k_1} \cdot (\pmb{\upalpha}^Q)^T\m{c}^{k_2}\right].
\end{equation}

Here $P = (\pi_1,\ldots,\pi_{p_1})$ ranges over set $\cC$ of all combinations of $p_1$ indices from $\mathcal{I} \defeq \{1,\ldots,p\}$, with $p = p_1 + p_2$. The complementary set is $Q \defeq \mathcal{I}\setminus P$. The number of terms in the sum is therefore ${p \choose p_1} = {p \choose p_2}$.

For each $P$, define the subvectors $\mt^P \defeq (t_{\pi_1}^i,\ldots,t_{\pi_{p_1}}^i)$ and $\mt^Q$ from $\mt^i \subseteq \mt$. Let $k_1$ and $k_2$ be such that $t_i \in [\tau_{k_\ell}^\ell,\tau_{k_\ell+1}^\ell)$ for $\ell=1,2$. Then define the local knot vectors $\mtau^{k_1}$, $\mtau^{k_2}$ and the local coefficient vectors $\m{c}^{k_1}$, $\m{c}^{k_2}$ as in the Oslo Algorithm. These data determine the discrete B-splines $\pmb{\upalpha}^P$ and $\pmb{\upalpha}^Q$ in \eqref{eq:diequation}.

An algorithmic implementation for computing $\m{b} = (b_i)_{i=1}^m$ is given in Algorithm \ref{alg:morken}.
\begin{algorithm}[h!]\fontsize{11.75}{11.75} \selectfont
\caption{M{\o}rken Algorithm - {\fontsize{10}{10} \selectfont$[p, \mt, \m{b}] = \texttt{mørken}(p_1, p_2, \mtau^1, \mtau^2, \m{c}^1, \m{c}^2)$}}\label{alg:morken}
\nonl\begin{minipage}{\textwidth}
\rule{0.95\textwidth}{0.025cm}

\textbf{Input:}\begin{itemize}\vspace{-0.2cm}\setlength{\itemindent}{0.45cm}
\item[$p_1, p_2\colon$] degrees of the spline factors\vspace{-0.05cm}
\item[$\pmb{\uptau}^1, \pmb{\uptau}^2\colon$] knot vectors of the spline factors\vspace{-0.05cm}
\item[$\m{c}^1, \m{c}^2\colon$] coefficient vectors of the spline factors\vspace{-0.05cm}
\end{itemize}
\noindent\textbf{Output:}\begin{itemize}\vspace{-0.2cm}\setlength{\itemindent}{-0.2cm}
\item[$p\colon$] degree of the spline product\vspace{-0.05cm}
\item[$\m{t}\colon$] knot vector of the spline product\vspace{-0.05cm}
\item[$\m{b}\colon$] coefficient vector of the spline product\vspace{-0.5cm}
\end{itemize}
\rule{0.95\textwidth}{0.025cm}\vspace{0.1cm}\end{minipage}

\linespread{1.35}\selectfont
\nl Set $p = p_1 + p_2$\;
Define $\m{t} = (t_j)_{j=1}^{m + p + 1}$ from $\pmb{\uptau}^1, \pmb{\uptau}^2$ as described in \eqref{eq:tmultiplicities}\;
Initialize $\m{b} = \pmb{0} \in \RR^m$\;
Let $\mathcal{C}$ be the set of combinations of $p_1$ indices from $\mathcal{I} \defeq \{1, \ldots, p\}$\;
\nonl\For{$i = 1, \ldots, m$}{
 Let $\m{t}^i \defeq (t_j)_{j=i+1}^{i+p}$\;
If $t_i \in [\tau_{k_\ell}^\ell, \tau_{k_\ell+1}^\ell)$, let $\pmb{\uptau}^{k_\ell} \defeq (\tau_j^\ell)_{j = k_\ell + 1 - p_\ell}^{k_\ell + p_\ell} \subseteq \pmb{\uptau}^\ell$, for $\ell = 1, 2$\;
Let $\m{c}^{k_\ell} \defeq (c_j^\ell)_{j=k_\ell - p_\ell}^{k_\ell} \subseteq \m{c}^\ell$, for $\ell = 1, 2$\;
\nonl\For{all $P \defeq (\pi_1, \ldots, \pi_{p_1}) \in \mathcal{C}$}{
Let $\m{t}^P \defeq (t_{i + \pi_k})_{k=1}^{p_1} \subseteq \m{t}^i$\;
Let $\m{t}^Q \defeq \m{t}^i \setminus \m{t}^P$\;
Set $b_i^P = \texttt{deboor}(p_1, \pmb{\uptau}^{k_1}, \m{c}^{k_1}, \m{t}^P)$\; 
Set $b_i^Q = \texttt{deboor}(p_2, \pmb{\uptau}^{k_2}, \m{c}^{k_2}, \m{t}^Q)$\;
Update $b_i \leftarrow b_i + b_i^P \cdot b_i^Q$\;
}}
Update $\m{b} \leftarrow \m{b} / {p \choose p_1}$\;
\end{algorithm}
The Oslo Algorithm appears explicitly in lines 5--7 and 10--11. Each de Boor call has cost $\frac{3}{2}p_\ell^2$, $\ell=1,2$. Hence, each update of $b_i$ costs about $\frac{3}{2}(p_1^2 + p_2^2)$. The total cost is
$$
\frac{3}{2}{p \choose p_1} m (p_1^2 + p_2^2).
$$

However, the computation of the coefficients $b_i$ can be parallelized. The outer loop over $i$ is independent and can be removed. The inner loop over combinations can also be parallelized, since all terms are independent. This only affects the order of summation.

Despite these parallelizations, the number of terms ${p \choose p_1}$ grows rapidly. The method is therefore practical only when $p_2$ is small relative to $p_1$. This limitation explains why the formula is rarely used in practice, where implicit methods such as collocation are preferred. In the next section we present an improved procedure that reduces the number of terms in the sum.

\begin{remark}
In the original formulation of \citet{morken}, the knot vectors $\mt^P$ and $\mt^Q$ are global and refine $\mtau^1$ and $\mtau^2$. Here we use local knot vectors. These may be non-open and span smaller intervals. They only need to refine $\mtau^1$ and $\mtau^2$ locally on their knot spans. This is valid due to the locality of the Oslo Algorithm and because openness is not required (see Remarks \ref{oss:notopen} and \ref{oss:localvec}).
\end{remark}

\subsection{The improved procedure}
The main drawback of Algorithm \ref{alg:morken} is the inner loop over all combinations of $p_1$ indices from a set of size $p$. This introduces a factor ${p \choose p_1}$ in the cost, which grows rapidly with $p$. Moreover, the algorithm does not exploit repeated knots in $\mt^i$. As a result, identical subvectors $\mt^P$ and $\mt^Q$ are recomputed multiple times.

To address this issue, we replace the loop over $\mathcal{C}$ with a loop over a set $\mathcal{T}^i$ of \emph{distinct} combinations derived from $\mt^i$. This set depends on $i$, unlike $\mathcal{C}$, which is global. We also introduce multiplicities. For each element of $\mathcal{T}^i$, we store the number of occurrences in the full set $\mathcal{C}$. These values form a set $\mathcal{M}^i$, paired with $\mathcal{T}^i$.

This leads to a factorization of the sum in \eqref{eq:diequation}. Identical terms are grouped together. Each distinct term corresponds to an element of $\mathcal{T}^i$ and is weighted by its multiplicity in $\mathcal{M}^i$. We denote by \texttt{knot\_combinations} the routine that computes $\mathcal{T}^i$ and $\mathcal{M}^i$ from $\mt^i$ and $p_1$. Its description is given in the next subsection.

The resulting procedure is shown in Algorithm \ref{alg:improved}. It computes $\m{b} = (b_i)_{i=1}^m$ more efficiently than Algorithm \ref{alg:morken}. The improvement is confirmed numerically in Section \ref{sec:numerics}.
\begin{algorithm}[h!]\fontsize{11.75}{11.75} \selectfont
\caption{ Improved M{\o}rken Algorithm - $[p, \mt, \m{b}] = \texttt{improved\_mørken}(p_1, p_2, \mtau^1, \mtau^2, \m{c}^1, \m{c}^2)$}\label{alg:improved}
\nonl\begin{minipage}{\textwidth}
\rule{0.95\textwidth}{0.025cm}

\textbf{Input:}\begin{itemize}\vspace{-0.2cm}\setlength{\itemindent}{0.45cm}
\item[$p_1, p_2\colon$] degrees of the spline factors\vspace{-0.05cm}
\item[$\pmb{\uptau}^1, \pmb{\uptau}^2\colon$] knot vectors of the spline factors\vspace{-0.05cm}
\item[$\m{c}^1, \m{c}^2\colon$] coefficient vectors of the spline factors\vspace{-0.05cm}
\end{itemize}
\noindent\textbf{Output:}\begin{itemize}\vspace{-0.2cm}\setlength{\itemindent}{-0.2cm}
\item[$p\colon$] degree of the spline product\vspace{-0.05cm}
\item[$\m{t}\colon$] knot vector of the spline product\vspace{-0.05cm}
\item[$\m{b}\colon$] coefficient vector of the spline product\vspace{-0.5cm}
\end{itemize}
\rule{0.95\textwidth}{0.025cm}\vspace{0.1cm}\end{minipage}

\linespread{1.35}\selectfont
\nl Set $p = p_1 + p_2$\;
Define $\m{t} = (t_j)_{j=1}^{m + p + 1}$ from $\pmb{\uptau}^1, \pmb{\uptau}^2$ as in \eqref{eq:tmultiplicities}\;
Initialize $\m{b} = \pmb{0} \in \RR^m$\;
\nonl\For{$i = 1, \ldots, m$}{
 Let $\m{t}^i \defeq (t_j)_{j=i+1}^{i+p}$\;
If $t_i \in [\tau_{k_\ell}^\ell, \tau_{k_\ell+1}^\ell)$, let $\pmb{\uptau}^{k_\ell} \defeq (\tau_j^\ell)_{j = k_\ell + 1 - p_\ell}^{k_\ell + p_\ell} \subseteq \pmb{\uptau}^\ell$, for $\ell = 1, 2$\;
Let $\m{c}^{k_\ell} \defeq (c_j^\ell)_{j=k_\ell - p_\ell}^{k_\ell} \subseteq \m{c}^\ell$, for $\ell = 1, 2$\;
Set $[\mathcal{T}^i, \mathcal{M}^i] = \texttt{knot\_combinations}(\mt^i, p_1)$\;
\nonl\For{all $(\m{t}^P, m^P)$ in $(\mathcal{T}^i, \cM^i)$}{
Let $\m{t}^Q \defeq \mt^i \setminus\mt^P$\;
Set $b_i^P = \texttt{deboor}(p_1, \pmb{\uptau}^{k_1}, \m{c}^{k_1}, \m{t}^P)$\; 
Set $b_i^Q = \texttt{deboor}(p_2, \pmb{\uptau}^{k_2}, \m{c}^{k_2}, \m{t}^Q)$\;
Update $b_i \leftarrow b_i + m^P \cdot b_i^P \cdot b_i^Q$\;
}}
Update $\m{b} \leftarrow \m{b} / {p \choose p_1}$\;
\end{algorithm}
The inner loop now iterates over $\mathcal{T}^i$, not over index combinations. We keep the notation $\mt^P$ for a generic element, to emphasize the link with Algorithm \ref{alg:morken}. Here, however, $P$ is only a label and does not range over $\mathcal{C}$.
Let
$$
\nu \defeq \sum_{i=1}^m |\cT^i|, \qquad \bar{\nu} \defeq \nu/m
$$
be the total number and mean number, respectively, of distinct combinations used to compute $\m{b}$. The total cost is then estimated as
$$
\frac{3}{2}\nu(p_1^2 + p_2^2) = \frac{3}{2}\bar{\nu} m (p_1^2 + p_2^2).
$$
This matches the cost of Algorithm \ref{alg:morken}, with ${p \choose p_1}$ replaced by $\bar{\nu}$. In Section \ref{sec:numerics} we compare these quantities and show the gain in efficiency.

The algorithm admits as well parallelization. The coefficients $b_i$ can be computed independently. For fixed $i$, each contribution associated with $\mt^P \in \mathcal{T}^i$ is also independent. Both loops can therefore be parallelized.

\begin{remark}
Both Algorithm \ref{alg:morken} and Algorithm \ref{alg:improved} require binomial coefficients. Their direct computation may overflow for large $p$. Two standard strategies can be used (see \citet{gamma}).

A first approach uses
$$
{n \choose k} = \prod_{j=1}^k \frac{n-k+j}{j},
$$
where we assumed $k = \min(k, n-k)$ without loss of generality. This form is stable when evaluated incrementally.
A second approach uses the Gamma function:
$$
{n \choose k} = \frac{\Gamma(n+1)}{\Gamma(k+1)\Gamma(n-k+1)}.
$$
In practice, one evaluates
$$
{n \choose k} = \exp\!\left(\log\Gamma(n+1) - \log\Gamma(k+1) - \log\Gamma(n-k+1)\right).
$$
This avoids overflow and replaces products by sums. It is the approach used in the experiments of Section \ref{sec:motivation}.
\end{remark}

\begin{remark}
A similar factorization idea appears in \citet{cohen}, applied to the formula of \citet{ueda}. The two approaches differ in representation. The formula of \citet{ueda} uses blossoms, while that of \citet{morken} uses discrete B-splines. 

The method of \citet{cohen} introduces a sliding window algorithm to evaluate blossoms efficiently. However, blossoming remains complex in practice (see \citet[Sections VII--VIII]{cohen}). 

Our approach proceeds differently. First, we express the formula of \citet{morken} in algorithmic form using a matrix-free procedure based on the Oslo Algorithm. Second, we factorize the expression to reduce the cost. This leads to a simpler and more practical method, easier to integrate into spline and isogeometric analysis codes.
\end{remark}

\subsection{Identification of distinct combinations and number of repetitions}\label{sec:distinctcomb}
For convenience, we write $\mt^i$ in terms of distinct breakpoints and multiplicities:
$$
\mt^i = \Theta^i \defeq (\theta_1^{m_1}, \ldots, \theta_s^{m_s}),
$$
where $\theta_j$ is the $j$th distinct knot and $m_j$ its multiplicity. Since $\mt^i$ has length $p$,
$$
\sum_{j=1}^s m_j = p.
$$

We construct all distinct combinations of $p_1$ knots from $\mt^i$ via a recursive procedure. The recursion is based on the number of breakpoints in $\Theta^i$. We start from the minimal case of at most two breakpoints. For longer sequences, we reduce the problem to a shorter sequence by removing one breakpoint.

Assume first that $\Theta^i$ has at most two breakpoints:
$$
\Theta^i = (\theta_1^{m_1}, \theta_2^{m_2}), \qquad m_1 + m_2 = p,
$$
with $m_2 = 0$ allowed. The multiplicities $\mu_1, \mu_2$ of a combination $(\theta_1^{\mu_1}, \theta_2^{\mu_2})$ of length $p_1$ satisfy
$$
\mu_1 + \mu_2 = p_1.
$$
The admissible range of $\mu_1$ is therefore
$$
\mu_1 \in \left[p_1 - \min\{p_1,m_2\},\; \min\{p_1,m_1\}\right].
$$
All distinct combinations are obtained by iterating over these values:

\begin{center}
\begin{minipage}{0.95\textwidth}
\begin{algorithm}[H]
\linespread{1.35}\selectfont
Let $m_{\max} \defeq \min\{p_1, m_1\}$ and $m_{\min} \defeq p_1 - \min\{p_1, m_2\}$\;
\nonl\For{$\mu_1 = m_{\max}, m_{\max}-1,\ldots, m_{\min}$}{
Set $\mu_2 \defeq p_1 - \mu_1$\;
$(\theta_1^{\mu_1}, \theta_2^{\mu_2})$ is a distinct combination\;
}
\end{algorithm}
\end{minipage}
\end{center}

If $\Theta^i$ has more than two breakpoints, we proceed recursively. The multiplicity $\mu_1$ of $\theta_1$ satisfies
$$
\mu_1 \in \left[p_1 - \min\!\left\{p_1, \sum_{\ell=2}^s m_\ell \right\},\; \min\{p_1,m_1\}\right].
$$
For each admissible $\mu_1$, we concatenate $\theta_1^{\mu_1}$ with all combinations of length $p_1-\mu_1$ from
$$
\widehat{\Theta}^i \defeq (\theta_2^{m_2}, \ldots, \theta_s^{m_s}).
$$
The same procedure is applied recursively to $\widehat{\Theta}^i$. The recursion terminates when only two breakpoints remain.

We now determine multiplicities. A distinct combination
$$
\mt^P = (\theta_1^{\mu_1^P}, \ldots, \theta_s^{\mu_s^P})
$$
is defined by the multiplicities $\mu_j^P$. For each breakpoint $\theta_j$, the number of ways to select $\mu_j^P$ occurrences among $m_j$ identical knots is ${m_j \choose \mu_j^P}$. Since selections are independent, the total multiplicity is
$$
m^P \defeq \prod_{j=1}^s {m_j \choose \mu_j^P}.
$$

This value counts how many times $\mt^P$ appears in the full set of combinations. It provides the weight used in Algorithm \ref{alg:improved} to group identical terms in \eqref{eq:diequation}. The full procedure is summarized in Algorithm \ref{alg:comb}.

\begin{algorithm}[h!]\fontsize{11.75}{11.75} \selectfont
\caption{{\fontsize{11.25}{11.25} \selectfont Distinct combinations and their multiplicities - $[\cT^i, \cM^i] = \texttt{knot\_combinations}(\m{t}^i, p_1)$}}\label{alg:comb}
\nonl\begin{minipage}{\textwidth}
\rule{0.95\textwidth}{0.025cm}

\textbf{Input:}\begin{itemize}\vspace{-0.2cm}\setlength{\itemindent}{-0.2cm}
\item[$\mt^i\colon$] $(t_j)_{j=i+1}^{i+p}\subseteq \mt$ local knot vector\vspace{-0.05cm}
\item[$p_1\colon$] spline factor degree\vspace{-0.05cm}
\end{itemize}
\noindent\textbf{Output:}\begin{itemize}\vspace{-0.2cm}\setlength{\itemindent}{-0.2cm}
\item[$\cT^i\colon$] collection of distinct combinations\vspace{-0.05cm}
\item[$\cM^i\colon$] multiplicities of the distinct combinations\vspace{-0.5cm}
\end{itemize}
\rule{0.95\textwidth}{0.025cm}\vspace{0.1cm}\end{minipage}

\linespread{1.35}\selectfont
\nl Let $\Theta^i \defeq (\theta_1^{m_1}, \ldots, \theta_s^{m_s})$ be the breakpoints in $\mt^i$ of multiplicities $m_1, \ldots, m_s$ with $m_1 > 0$ and $m_\ell \geq 0$ for $\ell = 2, \ldots, s$\;
Initialize $\cT^i \leftarrow \varnothing$ and $\cM^i \leftarrow \varnothing$\;
Set $m_{\max} \defeq \min\{p_1, m_1\}$\;
Set $m_{\min} \defeq p_1 - \min\!\left\{p_1, \sum_{\ell=2}^s m_\ell \right\}$\;
\nonl\uIf{$\Theta^i$ has at most two breakpoints, i.e., $s \leq 2$}{
\nonl    \For{$\mu_1 = m_{\max}, m_{\max}-1, \ldots, m_{\min}$}{
        Set $\mu_2 = p_1 - \mu_1$\;
        Insert $(\theta_1^{\mu_1}, \theta_2^{\mu_2})$ in $\cT^i$\vspace{0.1cm}\;
        Insert ${m_1 \choose \mu_1}{m_2 \choose \mu_2}$ in $\cM^i$\;
    }
}
\nonl\Else{
   \nonl \For{$\mu_1 = m_{\max}, m_{\max}-1, \ldots, m_{\min}$}{
        Let $\widehat{\mt}^i$ be $\mt^i$ with all the $m_1$ occurrences of $\theta_1$ removed\;
        Let $[\widehat{\cT}^i, \widehat{\cM}^i]
        = \texttt{knot\_combinations}(\widehat{\mt}^i, p_1 - \mu_1)$\;
        Concatenate $\theta_1^{\mu_1}$ to each element of $\widehat{\cT}^i$\;
        Insert the resulting knot vectors in $\cT^i$\;
        Insert in $\cM^i$ the rescaled repetitions ${m_1 \choose \mu_1}\cdot\widehat{\cM}^i$\;
    }
}
\end{algorithm}

\section{Computational aspects and comparisons on the considered numerical tests}\label{sec:numerics}
In this section we compare the number of terms required by the direct formula \eqref{eq:diequation} when implemented with Algorithm \ref{alg:morken} and with the improved Algorithm \ref{alg:improved}. 

For Algorithm~\ref{alg:morken}, the number of terms is fixed and equal to ${p \choose p_1}$. For Algorithm \ref{alg:improved}, it depends on the index $i$. We therefore consider the average number of terms, denoted by $\bar{\nu}$, over all indices $i$. The comparison is performed for the test cases of Section \ref{sec:motivation}.

The results show that Algorithm \ref{alg:improved} always requires fewer terms. In some cases the difference is dramatic. The binomial coefficient may reach values of order \texttt{1e29}, while $\bar{\nu}$ remains moderate. Even when the gap is smaller, the gain is still significant. This cost is incurred for each coefficient $b_i$. Hence, the advantage increases with the number of coefficients $m$.

The experiments also show that the type of factors (single B-splines or full splines) does not affect the number of terms per coefficient. It only changes the dimension of the product space $\SS_{\m{t}}^p$, i.e., the number of coefficients $m$ to compute. The comparison is therefore essentially unchanged.

For each test we report two figures: one comparing ${p \choose p_1}$ with $\bar{\nu}$ and one showing only the growth of $\bar{\nu}$.

\subsection{Product of a fixed B-spline, or spline, and polynomials of increasing degrees}\label{sec:splinepolyterms}
\begin{figure}
\centering
\subfloat[]{\includegraphics[width = 0.4\textwidth, page = 1]{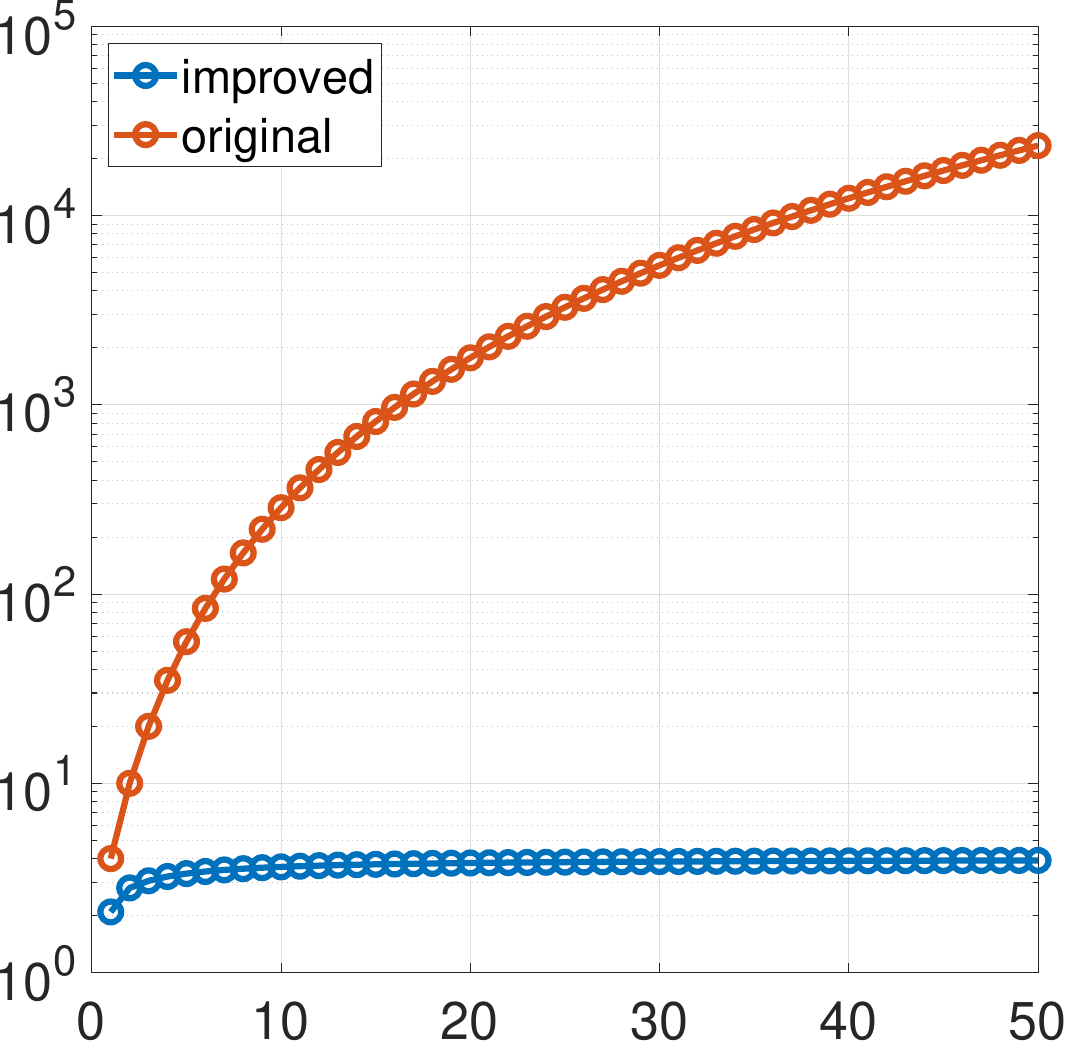}}\quad
\subfloat[]{\includegraphics[width = 0.4\textwidth, page = 2]{terms}}
\caption{Comparison of $\bar{\nu}$ and ${p \choose p_1}$ to compute one coefficient for the product of a fixed cubic B-spline and a polynomial of increasing degree. (a) comparison vs. polynomial degree. (b) $\bar{\nu}$ vs. polynomial degree.}
\label{fig:terms1}
\end{figure}

We first consider the case of Section \ref{sec:splinepoly}. A fixed cubic B-spline is multiplied by a polynomial of degree $p_2 \in \{1,\dots,50\}$. The results are shown in Figure \ref{fig:terms1} (a)--(b).
We also tested a general spline of degree $p_1 = 3$ with the same polynomial factor. The behavior is unchanged.

\subsection{The Galerkin isogeometric case with $p$- and $k$-refinements}
We next consider the product of two B-splines on the same knot vector $\pmb{\uptau}$, with equal degree $p_1 = p_2 \in \{3,\dots,50\}$. The knot multiplicities are adjusted to preserve global $C^2$ smoothness, as in $p$-refinement.

\begin{figure}
    \centering
    \subfloat[]{\includegraphics[width = 0.4\textwidth, page = 3]{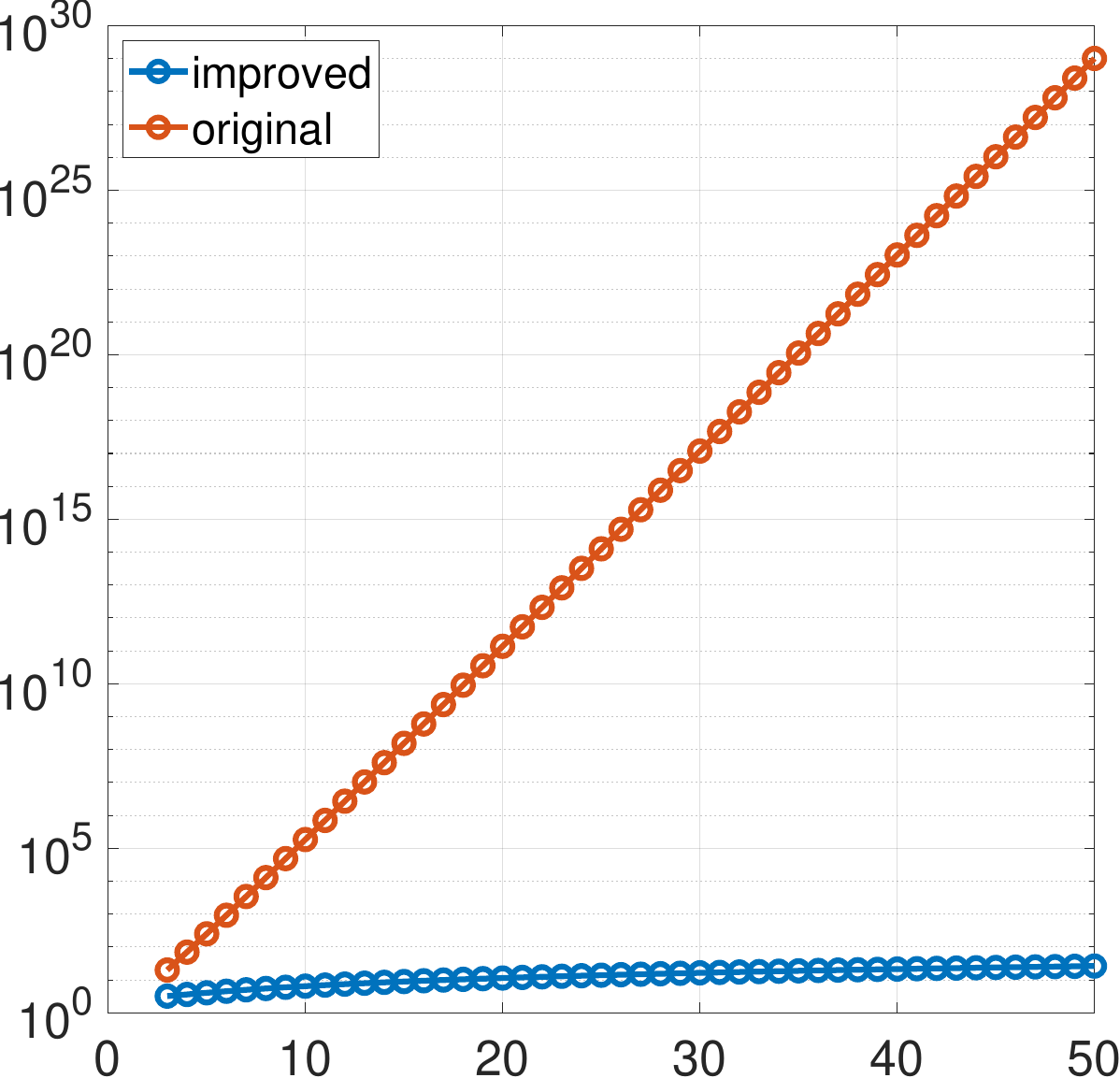}}\quad
    \subfloat[]{\includegraphics[width = 0.38\textwidth, page = 4]{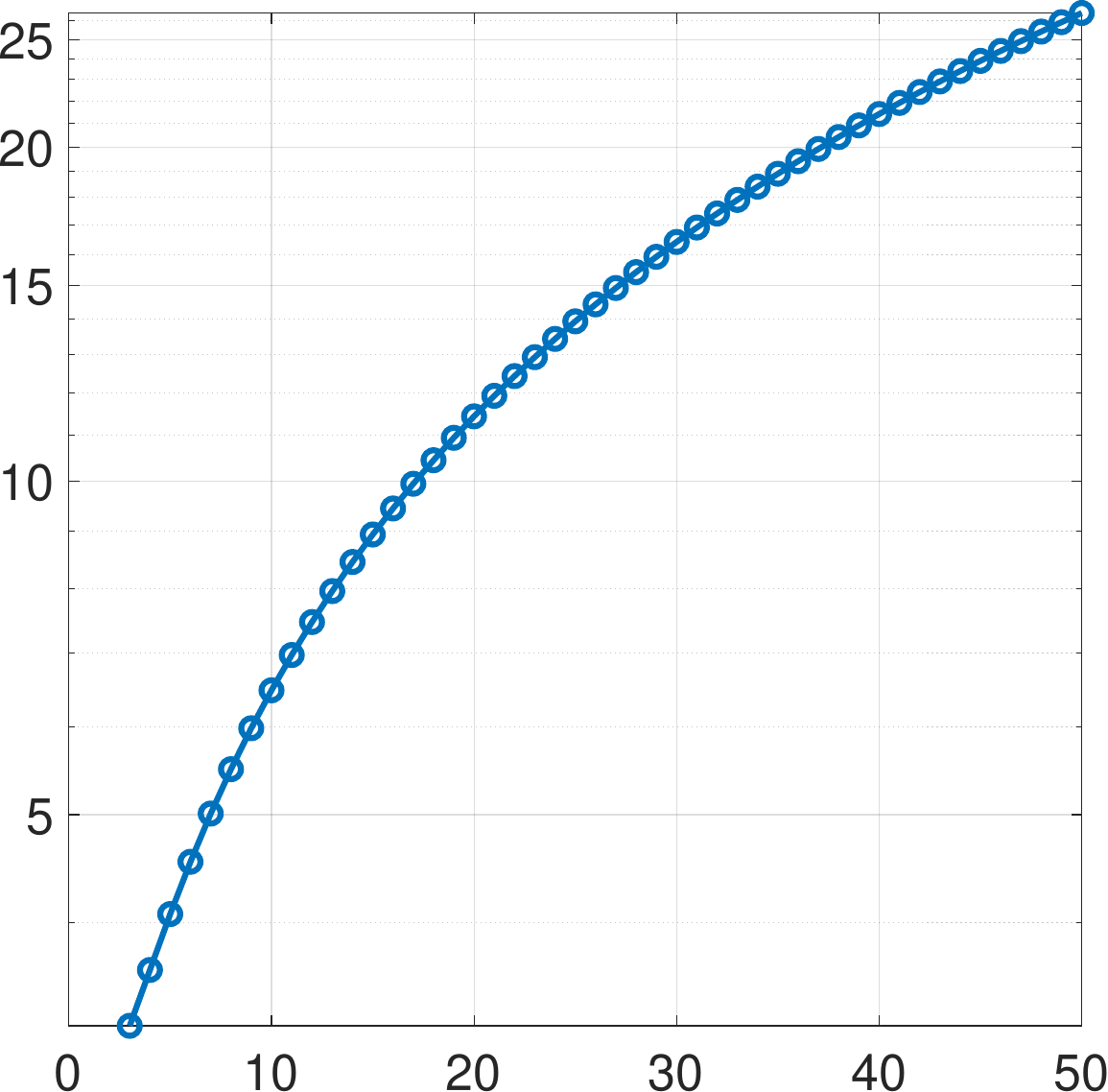}}
    \caption{Comparison of $\bar{\nu}$ (averaged over $i$ and $j$) and ${p \choose p_1}$ to compute one coefficient of $B_{i,\pmb{\uptau}} \cdot B_{j,\pmb{\uptau}}$ with fixed smoothness ($C^2$) and increasing degree. (a) comparison vs. polynomial degree. (b) $\bar{\nu}$  vs. polynomial degree.}
\label{fig:terms2}
\end{figure}

\begin{figure}
    \centering
    \subfloat[]{\includegraphics[width = 0.4\textwidth, page = 5]{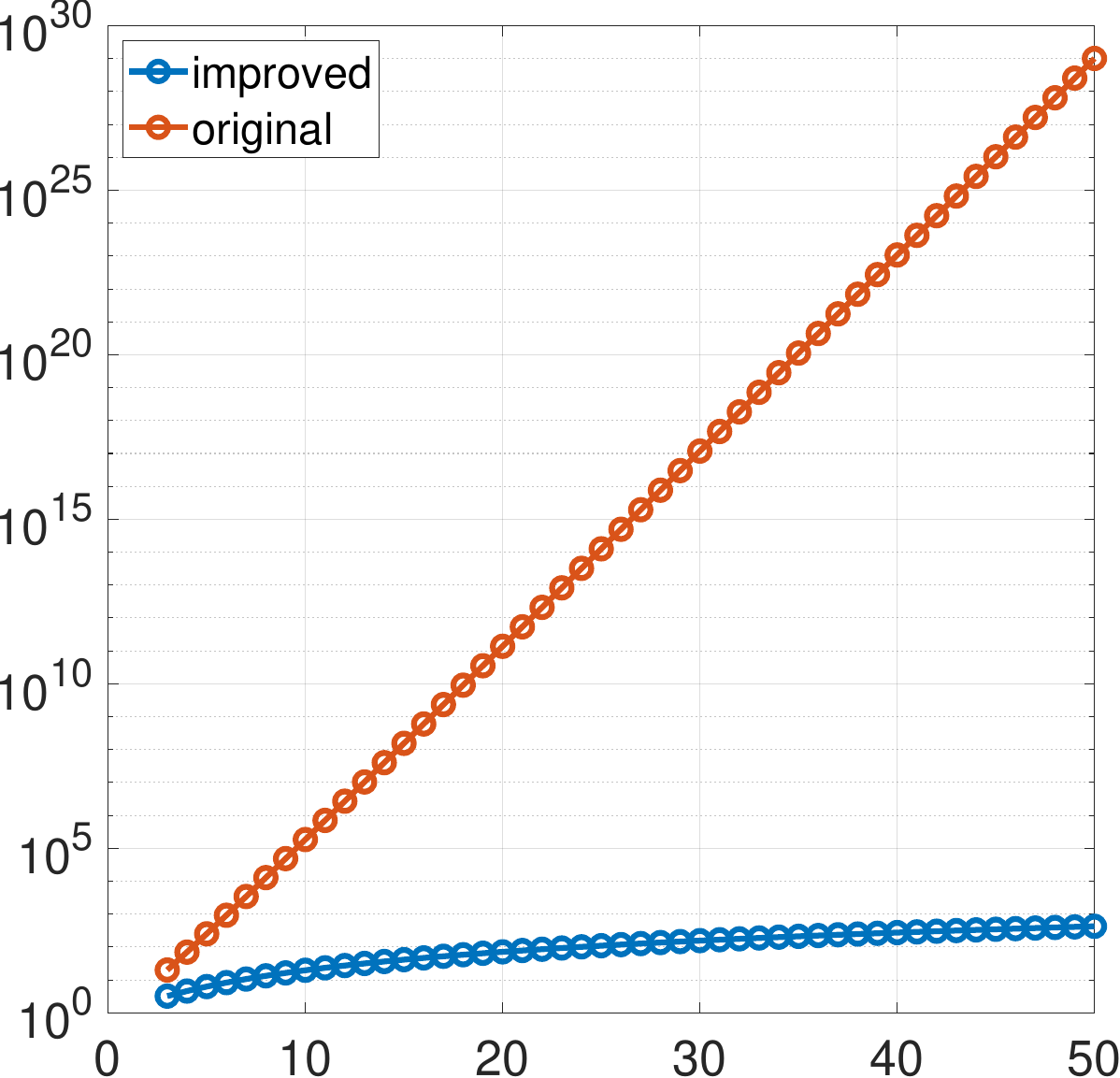}}\quad
    \subfloat[]{\includegraphics[width = 0.4\textwidth, page = 6]{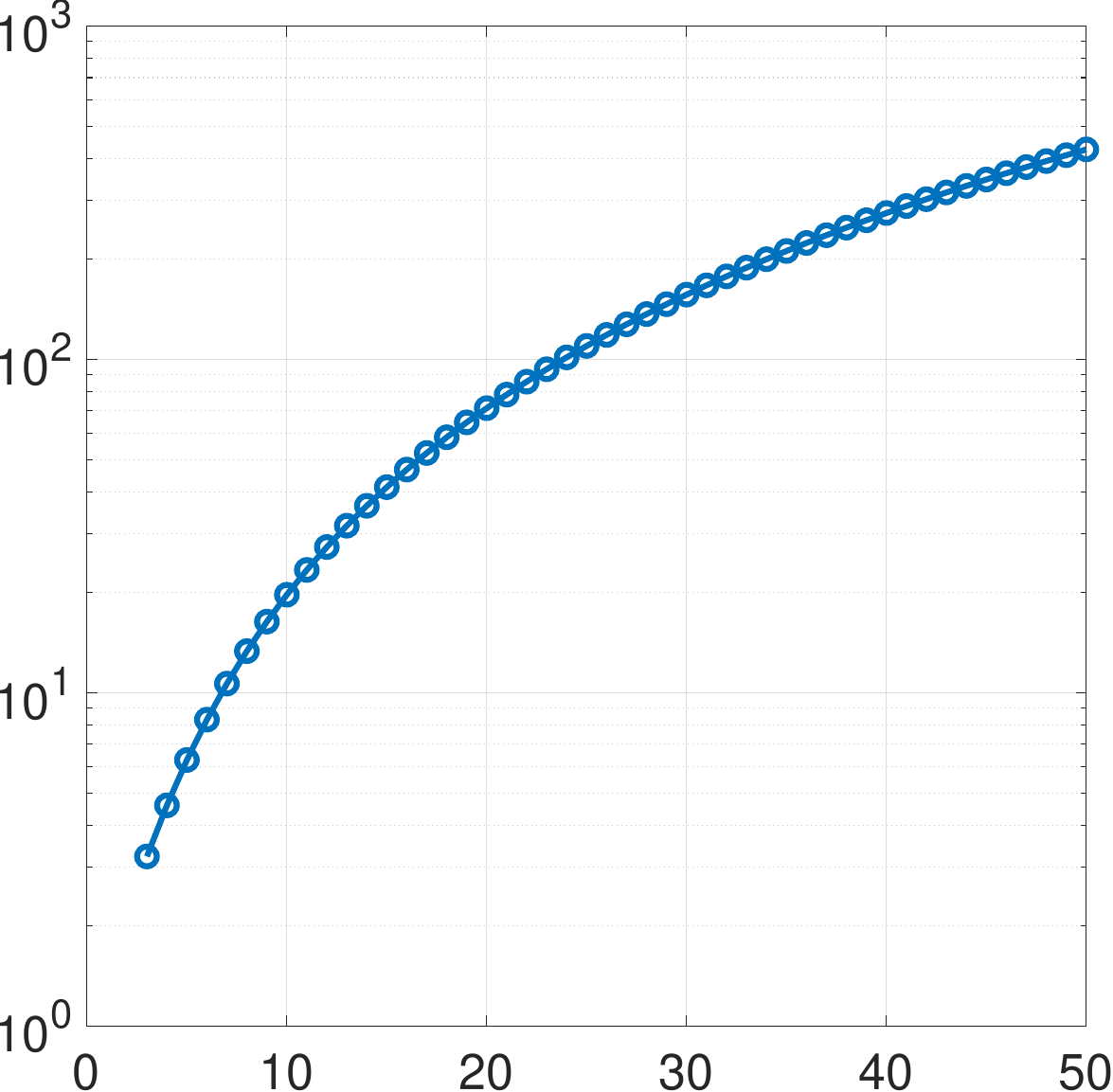}}
    \caption{Comparison of $\bar{\nu}$ (averaged over $i$ and $j$) and ${p \choose p_1}$ to compute one coefficient of $B_{i,\pmb{\uptau}} \cdot B_{j,\pmb{\uptau}}$ with maximal smoothness and increasing degree. (a) comparison vs. polynomial degree. (b) $\bar{\nu}$ vs. polynomial degree.}
\label{fig:terms21}
\end{figure}

For each degree, we fix $f = B_{i,\pmb{\uptau}}$ and consider all B-splines $B_{j,\pmb{\uptau}}$ with overlapping support, $\supp B_{i,\pmb{\uptau}} \cap \supp B_{j,\pmb{\uptau}} \neq \varnothing$. The average $\bar{\nu}$ is computed by averaging first over all such $j$, and then over $i=1,\dots,m$. The results are shown in Figure \ref{fig:terms2} (a)--(b).
The rapid growth of ${p \choose p_1}$ follows from the asymptotic behaviour of the central binomial coefficient, i.e., when $p_1 = \lceil p/2 \rceil$. By Stirling’s approximation, it behaves like $2^p/\sqrt{p}$ and grows exponentially with $p$.
We also consider the $k$-refinement case, where smoothness is maximal. The results are shown in Figure \ref{fig:terms21} (a)--(b).

\subsection{Product of splines on same knot vector and same degree}
\begin{figure}
\centering
\subfloat[]{\includegraphics[width = 0.4\textwidth, page = 7]{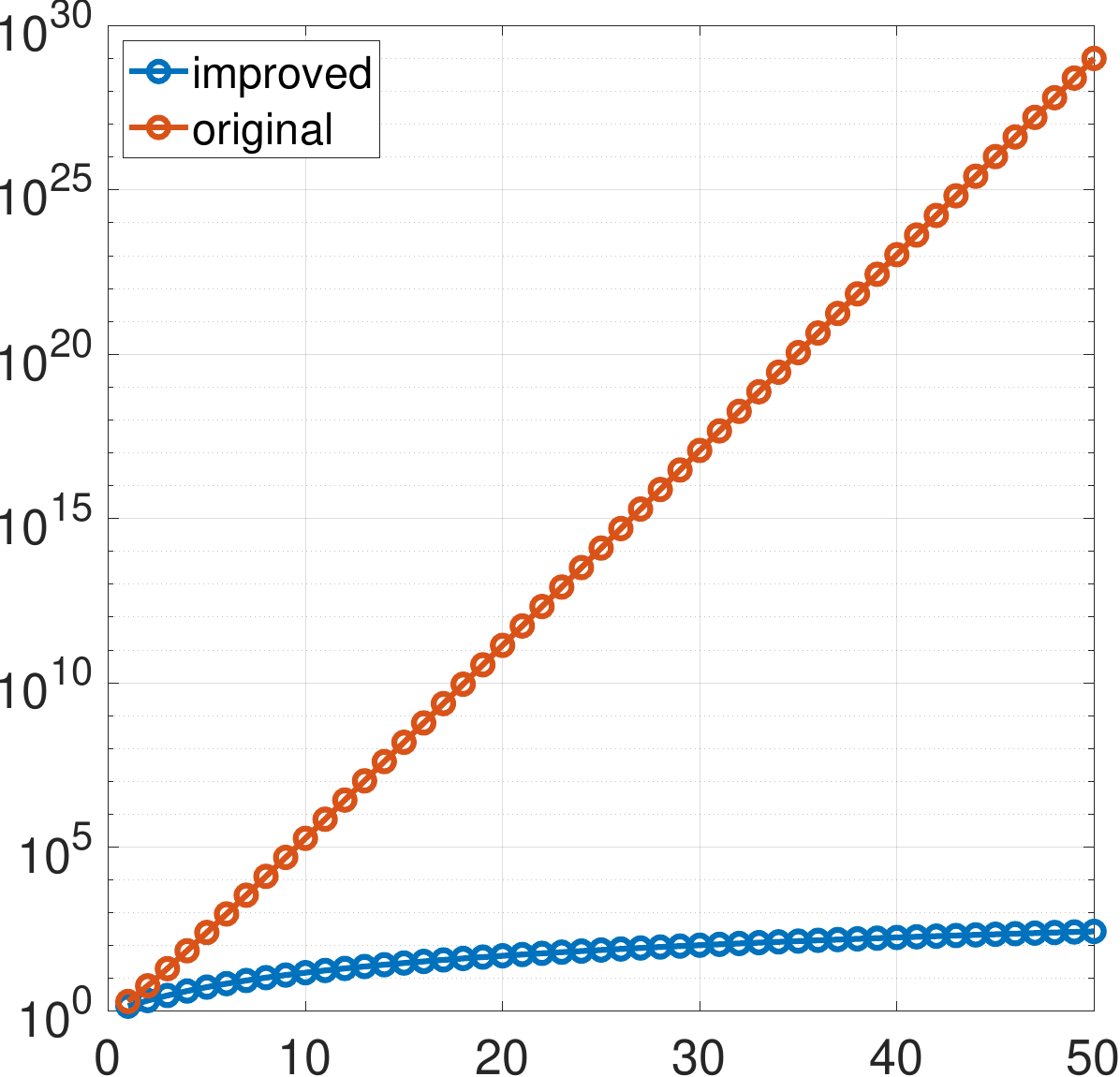}}\quad
\subfloat[]{\includegraphics[width = 0.4\textwidth, page = 8]{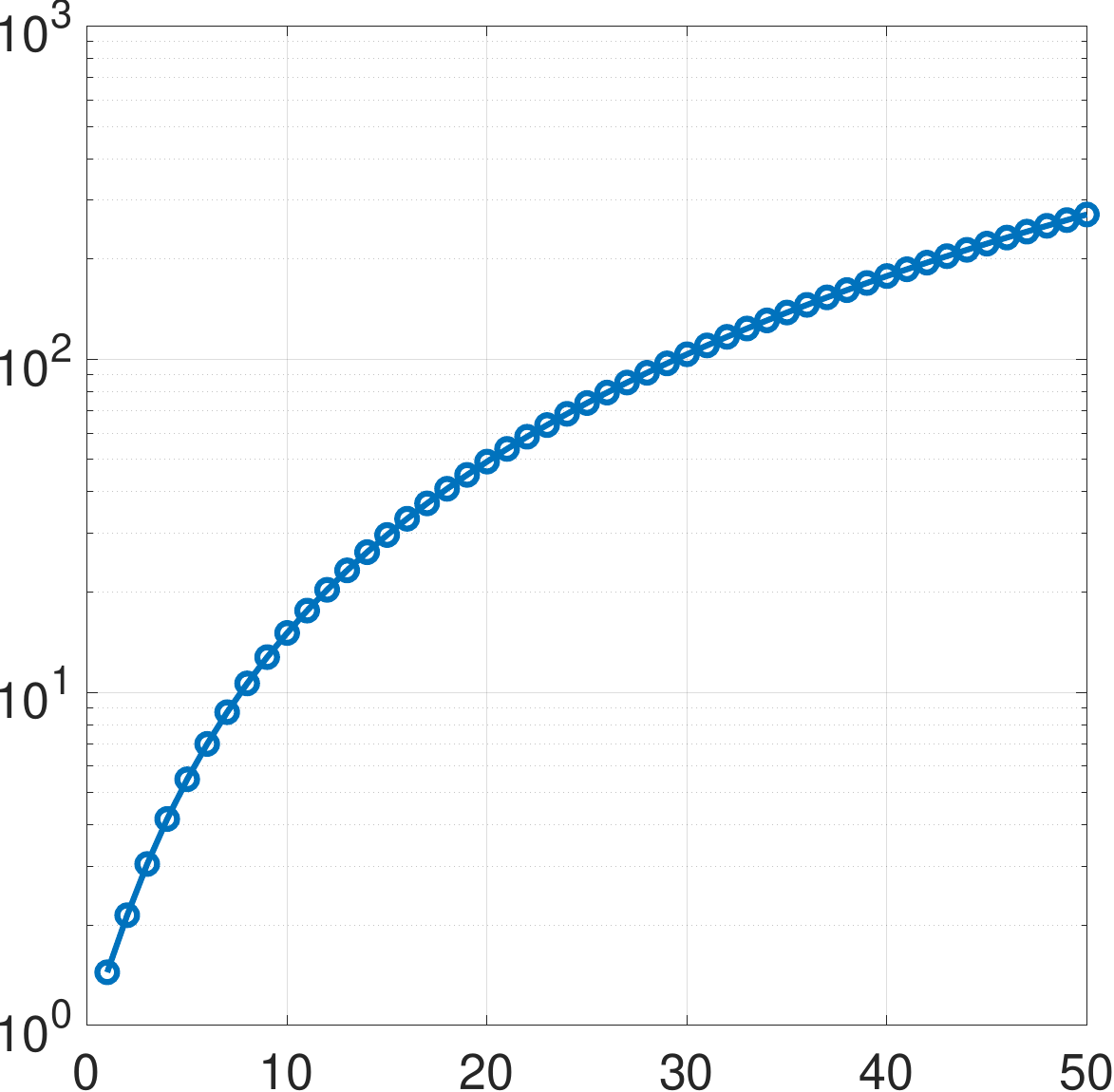}}
\caption{Comparison of $\bar{\nu}$ and ${p \choose p_1}$ to compute one coefficient for the product of splines on a common knot vector $\pmb{\uptau}$ with increasing degree. (a) comparison vs. polynomial degree. (b) $\bar{\nu}$ vs. polynomial degree.}
\label{fig:terms3}
\end{figure}

We consider the product of two splines on the same knot vector $\pmb{\uptau}$, with degree $p_1 = p_2 \in \{1,\dots,50\}$. The knot vector has five breakpoints in $[0,1]$. The coefficients of both splines are random in $[-1,1]$. The comparison is shown in Figure~\ref{fig:terms3} (a)--(b).

\subsection{Product of a fixed spline with splines on increasingly finer meshes}
\begin{figure}
\centering
\subfloat[]{\includegraphics[width = 0.39\textwidth, page = 9]{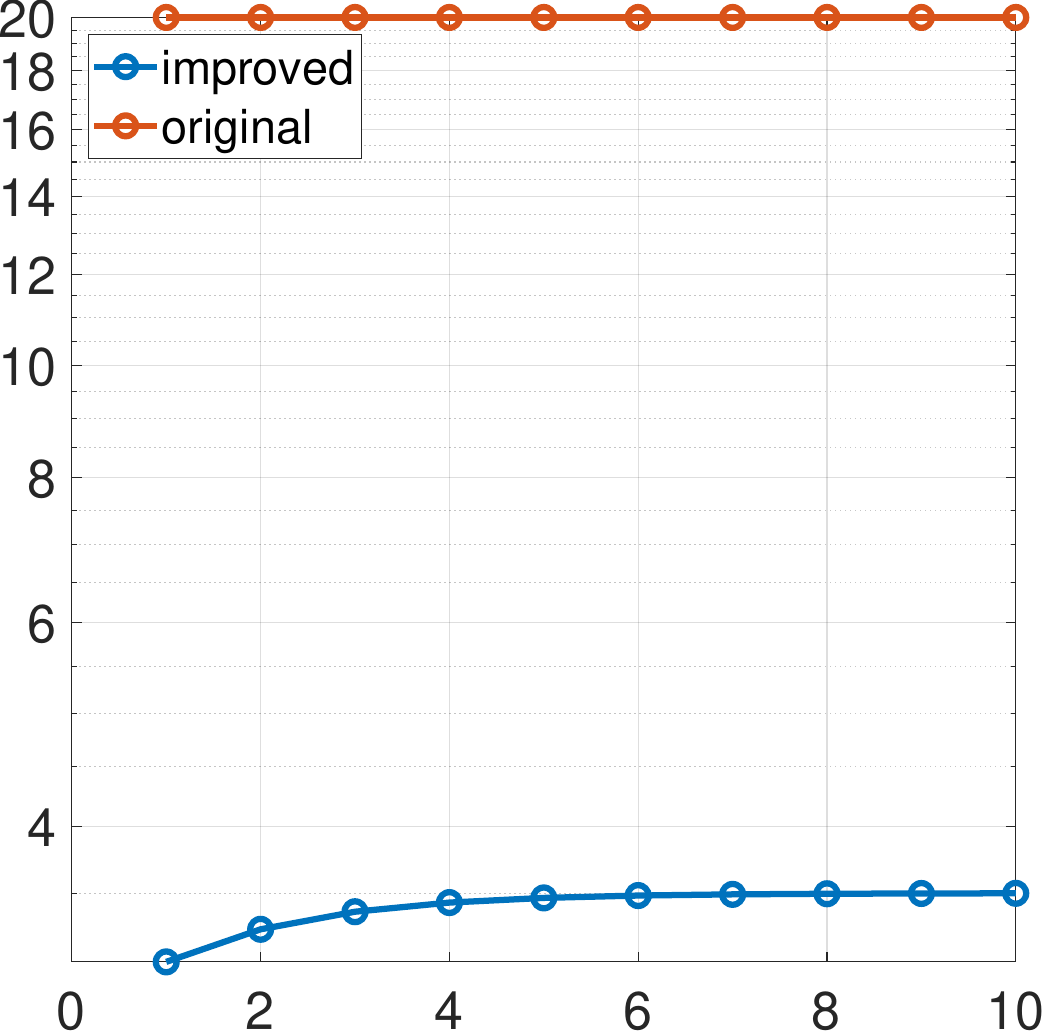}}\quad
\subfloat[]{\includegraphics[width = 0.41\textwidth, page = 10]{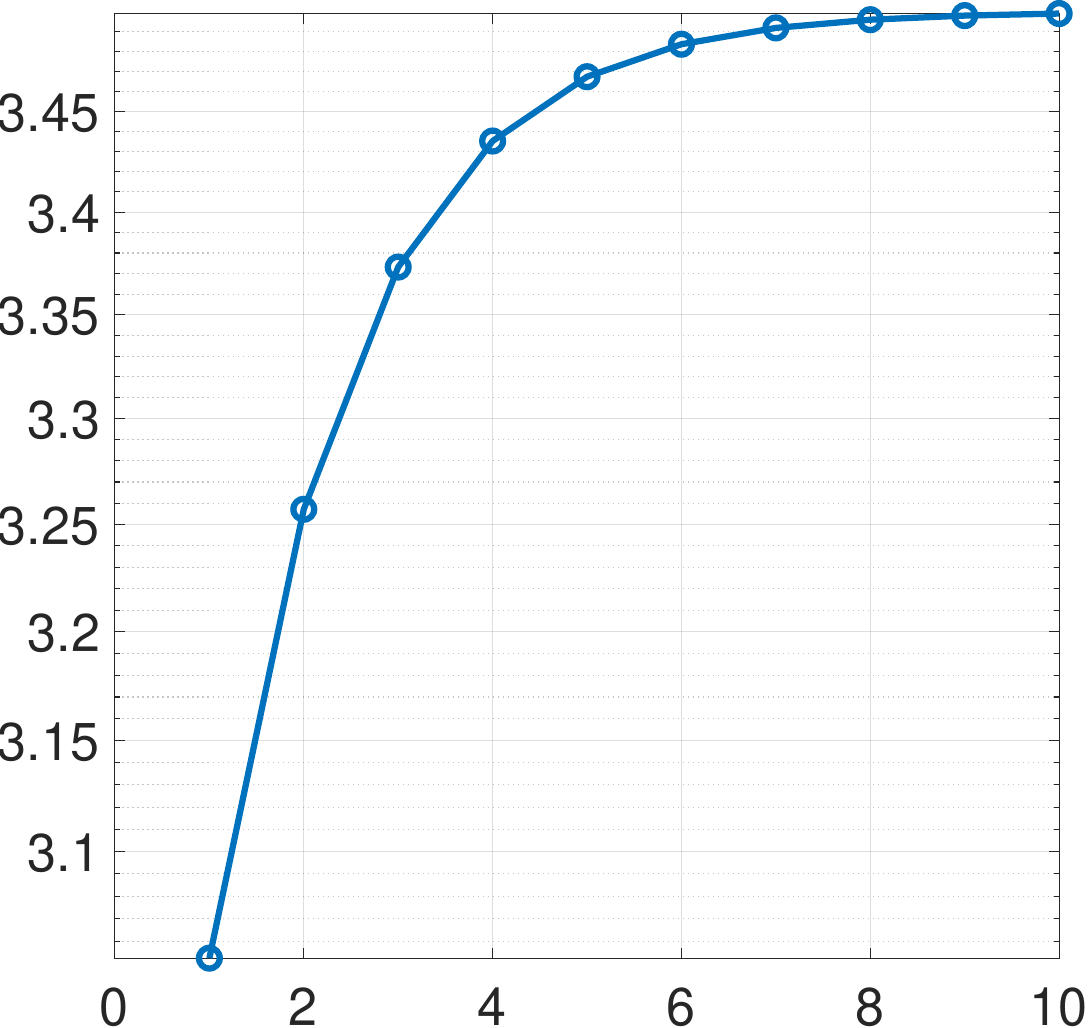}}
\caption{Comparison of $\bar{\nu}$ and ${p \choose p_1}$ to compute one coefficient for the product of cubic splines on nested meshes. The first spline is defined on a fixed knot vector $\mtau^1$ with $5$ uniform breakpoints, the second on $\mtau^2$ with $2^n+3$ breakpoints, $n=1,\ldots,10$. (a) comparison vs. $n$. (b) $\bar{\nu}$ vs. $n$.}
\label{fig:terms4}
\end{figure}

\begin{figure}
\centering
\subfloat[]{\includegraphics[width = 0.36\textwidth, page = 11]{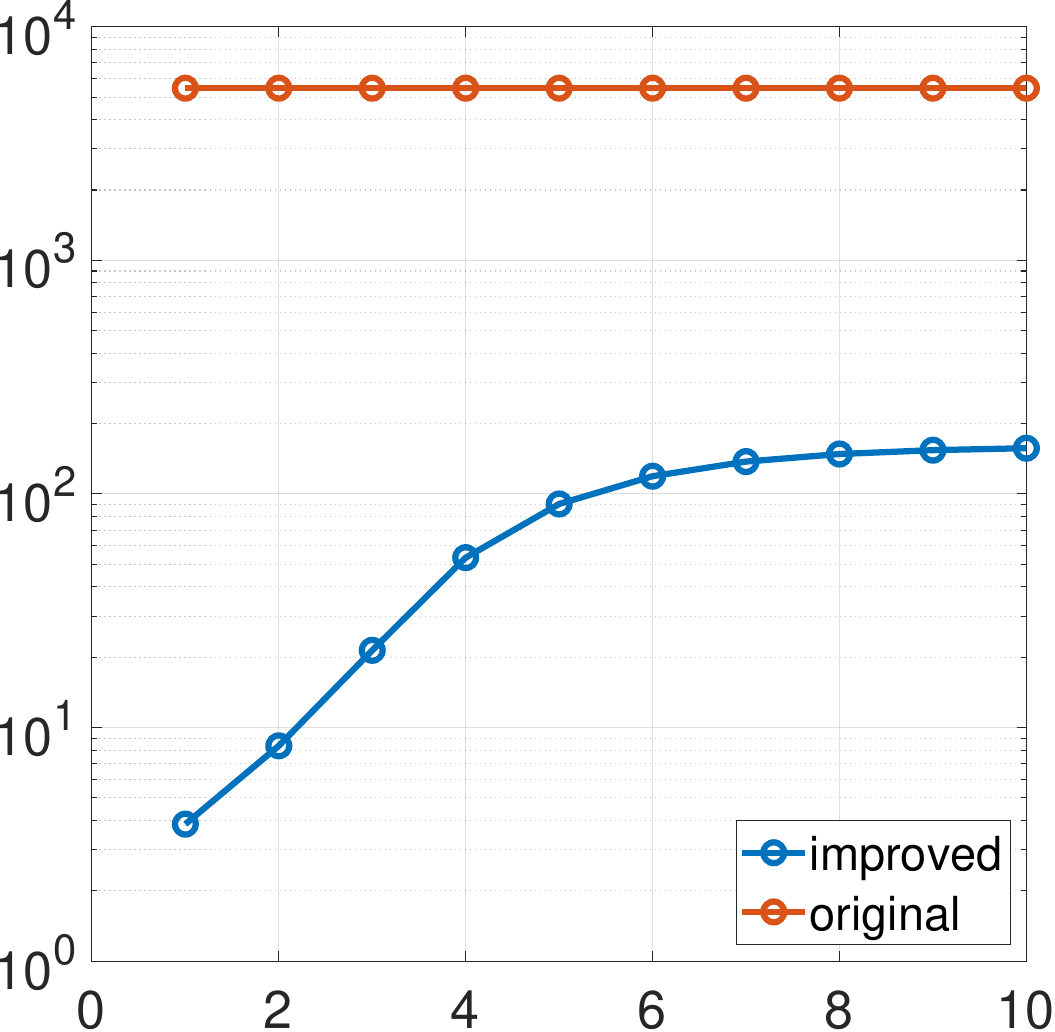}}\qquad
\subfloat[]{\includegraphics[width = 0.36\textwidth, page = 12]{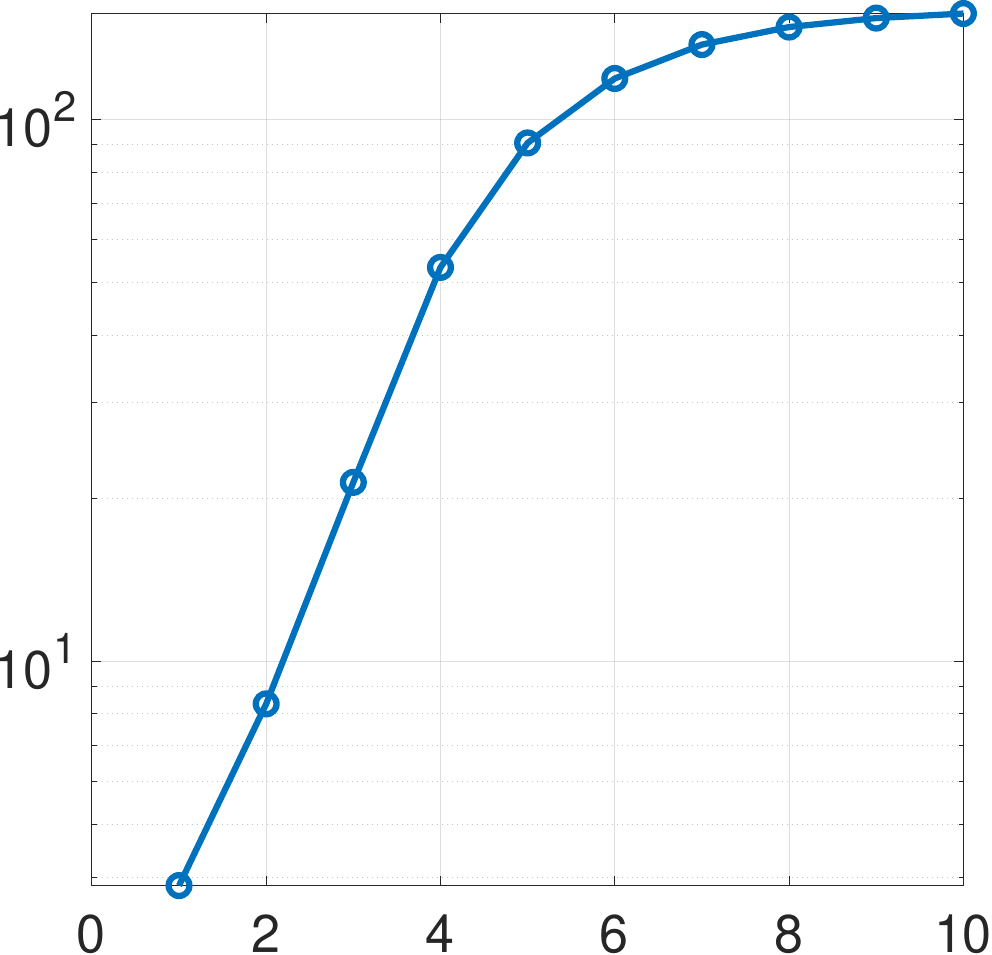}}
\caption{Comparison of $\bar{\nu}$ and ${p \choose p_1}$ to compute one coefficient for the product of a cubic spline and a spline of degree $30$ on nested meshes. The first spline is defined on a fixed knot vector $\mtau^1$ with $5$ uniform breakpoints, the second on $\mtau^2$ with $2^n+3$ breakpoints, $n=1,\ldots,10$. (a) comparison vs $n$. (b) $\bar{\nu}$ vs.\ $n$.}
\label{fig:terms41}
\end{figure}

Finally, we consider the product of a fixed spline of degree $p_1 = 3$ on a knot vector $\mtau^1$ with splines defined on successive refinements of $\mtau^1$.
Even for $p_2 = 3$, the improved procedure already yields a significant reduction in the number of terms, see Figure \ref{fig:terms4} (a)--(b). The advantage increases with the degree. For instance, for $p_2 = 30$, the gain is much more pronounced, see Figure \ref{fig:terms41} (a)--(b).

\section{Conclusion}\label{sec:conclusion}
In this work we have revisited the direct formula of \citet{morken} for representing spline products in the B-spline basis and shown that, when reformulated algorithmically, it yields a robust, accurate, and computationally efficient alternative to the implicit interpolation approach. By relying on the Oslo Algorithm, the computation is reduced to successive matrix–vector products with bidiagonal matrices. Then, by factoring out repeated contributions in the combinatorial sum, we achieve a substantial reduction in computational cost while preserving numerical stability.

Extensive numerical experiments confirm that the proposed algorithm remains efficient even for high polynomial degrees, where naive direct implementations of the direct formula become impractical and implicit methods suffer from severe ill-conditioning. In representative test cases, the procedure is often superior in accuracy by several orders of magnitude compared to the interpolation approach. 

Several directions remain open for future research. First, while the algorithm exhibits an intrinsic parallel structure, both across spline coefficients and across distinct knot combinations, a systematic investigation of parallel and high-performance implementations, including shared-memory architectures, is of clear interest for scalability. Second, a deeper theoretical analysis of the computational complexity of the improved procedure, including sharp bounds and average-case estimates, would further clarify its advantages over blossoming-based and implicit approaches. 
Third, it would be of interest to establish direct formulations for other fundamental spline operations, such as the composition of splines, and extend the direct-versus-implicit comparison carried out in this work.
Finally, adaptive variants of the algorithm, in which negligible contributions to the computation of the product B-spline coefficients are automatically dismissed, represent a possible direction for further reducing computational effort in large-scale applications.

Overall, the results of this work suggest that direct spline products, long regarded as theoretically elegant but practically cumbersome, can in fact be implemented in a simple, efficient, and robust manner, making them a viable and attractive tool for a wide range of spline-based numerical methods.

\section*{Acknowledgments}
This work has been supported by the MUR Excellence Department Project MatMod@TOV (CUP E83C23000330006) awarded to the Department of Mathematics of the University of Rome Tor Vergata.
The authors are members of \textit{Gruppo Nazionale per il Calcolo Scientifico} of the Italian \textit{Istituto Nazionale di Alta Matematica “Francesco Severi”} (GNCS-INdAM). The INdAM support through GNCS-2026 annual project (CUP E53C25002010001) is gratefully acknowledged.




 \bibliographystyle{elsarticle-harv} 
 \bibliography{biblio}

@article {morken,
    AUTHOR = {M{\o}rken, K.},
     TITLE = {Some identities for products and degree raising of splines},
   JOURNAL = {Constr. Approx.},
  FJOURNAL = {Constructive Approximation. An International Journal for
              Approximations and Expansions},
    VOLUME = {7},
      YEAR = {1991},
    NUMBER = {2},
     PAGES = {195--208},
      ISSN = {0176-4276,1432-0940},
   MRCLASS = {41A15 (65D07)},
  MRNUMBER = {1101062},
MRREVIEWER = {Ewald\ Quak},
       DOI = {10.1007/BF01888153},
       URL = {https://doi.org/10.1007/BF01888153},
}

@incollection {manni1,
    AUTHOR = {Lyche, Tom and Manni, Carla and Speleers, Hendrik},
     TITLE = {Foundations of spline theory: {B}-splines, spline
              approximation, and hierarchical refinement},
 BOOKTITLE = {Splines and {PDE}s: from approximation theory to numerical
              linear algebra},
    SERIES = {Lecture Notes in Math.},
    VOLUME = {2219},
     PAGES = {1--76},
 PUBLISHER = {Springer, Cham},
      YEAR = {2018},
      ISBN = {978-3-319-94910-9; 978-3-319-94911-6},
   MRCLASS = {65D07 (41A15)},
  MRNUMBER = {3839186},
MRREVIEWER = {Peter\ Alfeld},
}

@incollection {manni2,
    AUTHOR = {Manni, Carla and Speleers, Hendrik},
     TITLE = {Standard and non-standard {CAGD} tools for isogeometric
              analysis: a tutorial},
 BOOKTITLE = {Isogeometric analysis: a new paradigm in the numerical
              approximation of {PDE}s},
    SERIES = {Lecture Notes in Math.},
    VOLUME = {2161},
     PAGES = {1--69},
 PUBLISHER = {Springer, [Cham]},
      YEAR = {2016},
      ISBN = {978-3-319-42308-1; 978-3-319-42309-8},
   MRCLASS = {65D07 (65D17 65N30)},
  MRNUMBER = {3586483},
}

@book {deboor,
    AUTHOR = {de Boor, Carl},
     TITLE = {A practical guide to splines},
    SERIES = {Applied Mathematical Sciences},
    VOLUME = {27},
 PUBLISHER = {Springer-Verlag, New York-Berlin},
      YEAR = {1978},
     PAGES = {xxiv+392},
      ISBN = {0-387-90356-9},
   MRCLASS = {65D07 (41A15 65-01)},
  MRNUMBER = {507062},
MRREVIEWER = {Gerhard\ Merz},
}

@book {schumaker,
    AUTHOR = {Schumaker, Larry L.},
     TITLE = {Spline functions: basic theory},
    SERIES = {Cambridge Mathematical Library},
   EDITION = {Third},
 PUBLISHER = {Cambridge University Press, Cambridge},
      YEAR = {2007},
     PAGES = {xvi+582},
      ISBN = {978-0-521-70512-7},
   MRCLASS = {41-02 (41A15 65D07)},
  MRNUMBER = {2348176},
       DOI = {10.1017/CBO9780511618994},
       URL = {https://doi.org/10.1017/CBO9780511618994},
}

@article{piegl,
  title={Symbolic operators for {NURBS}},
  author={Piegl, Les and Tiller, Wayne},
  journal={Computer-Aided Design},
  volume={29},
  number={5},
  pages={361--368},
  year={1997},
  publisher={Elsevier}
}

@incollection {ueda,
    AUTHOR = {Ueda, Kenji},
     TITLE = {Multiplication as a general operation for splines},
 BOOKTITLE = {Curves and surfaces in geometric design
              ({C}hamonix-{M}ont-{B}lanc, 1993)},
     PAGES = {475--482},
 PUBLISHER = {A K Peters, Wellesley, MA},
      YEAR = {1994},
      ISBN = {1-56881-039-3},
   MRCLASS = {65D07 (65Y25)},
  MRNUMBER = {1302229},
}

@article {lee,
    AUTHOR = {Lee, E. T. Y.},
     TITLE = {Computing a chain of blossoms, with application to products of
              splines},
   JOURNAL = {Comput. Aided Geom. Design},
  FJOURNAL = {Computer Aided Geometric Design},
    VOLUME = {11},
      YEAR = {1994},
    NUMBER = {6},
     PAGES = {597--620},
      ISSN = {0167-8396,1879-2332},
   MRCLASS = {65D17 (65D07)},
  MRNUMBER = {1305909},
MRREVIEWER = {Richard\ D.\ Fuhr},
       DOI = {10.1016/0167-8396(94)90054-X},
       URL = {https://doi.org/10.1016/0167-8396(94)90054-X},
}

@article{cohen,
  title={An algorithm for direct multiplication of {B}-splines},
  author={Chen, Xianming and Riesenfeld, Richard F and Cohen, Elaine},
  journal={{IEEE} transactions on automation science and engineering},
  volume={6},
  number={3},
  pages={433--442},
  year={2009},
  publisher={{IEEE}}
}

@article {rida,
    AUTHOR = {Farouki, Rida T.},
     TITLE = {The {B}ernstein polynomial basis: a centennial retrospective},
   JOURNAL = {Comput. Aided Geom. Design},
  FJOURNAL = {Computer Aided Geometric Design},
    VOLUME = {29},
      YEAR = {2012},
    NUMBER = {6},
     PAGES = {379--419},
      ISSN = {0167-8396,1879-2332},
   MRCLASS = {65-03 (01A60 41A10 65D05 65D17)},
  MRNUMBER = {2921860},
       DOI = {10.1016/j.cagd.2012.03.001},
       URL = {https://doi.org/10.1016/j.cagd.2012.03.001},
}

@article{oslo1,
  title={Discrete {B}-splines and subdivision techniques in computer-aided geometric design and computer graphics},
  author={Cohen, Elaine and Lyche, Tom and Riesenfeld, Richard},
  journal={Computer graphics and image processing},
  volume={14},
  number={2},
  pages={87--111},
  year={1980},
  publisher={Elsevier}
}

@article {oslo2,
    AUTHOR = {Lyche, T. and M{\o}rken, K.},
     TITLE = {Making the {O}slo algorithm more efficient},
   JOURNAL = {SIAM J. Numer. Anal.},
  FJOURNAL = {SIAM Journal on Numerical Analysis},
    VOLUME = {23},
      YEAR = {1986},
    NUMBER = {3},
     PAGES = {663--675},
      ISSN = {0036-1429},
   MRCLASS = {65D07 (41A15)},
  MRNUMBER = {842650},
MRREVIEWER = {E.\ L.\ Albasiny},
       DOI = {10.1137/0723042},
       URL = {https://doi.org/10.1137/0723042},
}

@article{oslo3,
  title={Note on the {O}slo algorithm},
  author={Lyche, Tom},
  journal={Computer-aided design},
  volume={20},
  number={6},
  pages={353--355},
  year={1988},
  publisher={Elsevier}
}

@book{oslo4,
  title = {Spline {M}ethods {D}raft},
  author = {Lyche, Tom and M{\o}rken, Knut},
  year = {2018}
}

@article {prefinement1,
    AUTHOR = {Hughes, T. J. R. and Cottrell, J. A. and Bazilevs, Y.},
     TITLE = {Isogeometric analysis: {CAD}, finite elements, {NURBS}, exact
              geometry and mesh refinement},
   JOURNAL = {Comput. Methods Appl. Mech. Engrg.},
  FJOURNAL = {Computer Methods in Applied Mechanics and Engineering},
    VOLUME = {194},
      YEAR = {2005},
    NUMBER = {39-41},
     PAGES = {4135--4195},
      ISSN = {0045-7825,1879-2138},
   MRCLASS = {65D17 (65N30 74S05)},
  MRNUMBER = {2152382},
       DOI = {10.1016/j.cma.2004.10.008},
       URL = {https://doi.org/10.1016/j.cma.2004.10.008},
}

@article {prefinement2,
    AUTHOR = {Cottrell, J. A. and Hughes, T. J. R. and Reali, A.},
     TITLE = {Studies of refinement and continuity in isogeometric structural analysis},
   JOURNAL = {Comput. Methods Appl. Mech. Engrg.},
  FJOURNAL = {Computer Methods in Applied Mechanics and Engineering},
    VOLUME = {196},
      YEAR = {2007},
    NUMBER = {41},
     PAGES = {4160--4183},
      ISSN = {0045-7825},
       DOI = {10.1016/j.cma.2007.04.007},
}

@article {prefinement3,
    AUTHOR = {Beir\~ao da Veiga, L. and Buffa, A. and Rivas, J. and
              Sangalli, G.},
     TITLE = {Some estimates for {$h$}-{$p$}-{$k$}-refinement in
              isogeometric analysis},
   JOURNAL = {Numer. Math.},
  FJOURNAL = {Numerische Mathematik},
    VOLUME = {118},
      YEAR = {2011},
    NUMBER = {2},
     PAGES = {271--305},
      ISSN = {0029-599X,0945-3245},
   MRCLASS = {65D17 (41A15 41A20 41A63 65D07)},
  MRNUMBER = {2800710},
MRREVIEWER = {Manfred\ Tasche},
       DOI = {10.1007/s00211-010-0338-z},
       URL = {https://doi.org/10.1007/s00211-010-0338-z},
}

@article {prefinement4,
    AUTHOR = {Sangalli, G. and Tani, M.},
     TITLE = {Matrix-free weighted quadrature for a computationally
              efficient isogeometric {$k$}-method},
   JOURNAL = {Comput. Methods Appl. Mech. Engrg.},
  FJOURNAL = {Computer Methods in Applied Mechanics and Engineering},
    VOLUME = {338},
      YEAR = {2018},
     PAGES = {117--133},
      ISSN = {0045-7825,1879-2138},
   MRCLASS = {65N30 (65D07 65N22)},
  MRNUMBER = {3811621},
       DOI = {10.1016/j.cma.2018.04.029},
       URL = {https://doi.org/10.1016/j.cma.2018.04.029},
}

@article {prefinement5,
    AUTHOR = {Sande, Espen and Manni, Carla and Speleers, Hendrik},
     TITLE = {Sharp error estimates for spline approximation: explicit
              constants, {$n$}-widths, and eigenfunction convergence},
   JOURNAL = {Math. Models Methods Appl. Sci.},
  FJOURNAL = {Mathematical Models and Methods in Applied Sciences},
    VOLUME = {29},
      YEAR = {2019},
    NUMBER = {6},
     PAGES = {1175--1205},
      ISSN = {0218-2025,1793-6314},
   MRCLASS = {65D07 (41A15 41A27 41A44)},
  MRNUMBER = {3963638},
MRREVIEWER = {Dietrich\ Braess},
       DOI = {10.1142/S0218202519500192},
       URL = {https://doi.org/10.1142/S0218202519500192},
}

@article {prefinement6,
    AUTHOR = {Manni, Carla and Sande, Espen and Speleers, Hendrik},
     TITLE = {Outlier-free spline spaces for isogeometric discretizations of
              biharmonic and polyharmonic eigenvalue problems},
   JOURNAL = {Comput. Methods Appl. Mech. Engrg.},
  FJOURNAL = {Computer Methods in Applied Mechanics and Engineering},
    VOLUME = {417},
      YEAR = {2023},
     PAGES = {Paper No. 116314, 33},
      ISSN = {0045-7825,1879-2138},
   MRCLASS = {65N25 (65N30)},
  MRNUMBER = {4668245},
MRREVIEWER = {Jiansong\ Deng},
       DOI = {10.1016/j.cma.2023.116314},
       URL = {https://doi.org/10.1016/j.cma.2023.116314},
}

@book{gamma,
  author    = {William H. Press and Saul A. Teukolsky and William T. Vetterling and Brian P. Flannery},
  title     = {{N}umerical {R}ecipes: {T}he {A}rt of {S}cientific {C}omputing},
  publisher = {Cambridge University Press},
  address   = {Cambridge, UK},
  edition   = {3rd},
  year      = {2007},
  isbn      = {978-0-521-88068-8}
}

@article {stokes,
    AUTHOR = {Bracco, Cesare and Patrizi, Francesco and Sestini, Alessandra},
     TITLE = {A smoothly varying quadrature approach for 3{D} {I}g{A}-{BEM}
              discretizations: application to {S}tokes flow simulations},
   JOURNAL = {Comput. Methods Appl. Mech. Engrg.},
  FJOURNAL = {Computer Methods in Applied Mechanics and Engineering},
    VOLUME = {452},
      YEAR = {2026},
     PAGES = {Paper No. 118773, 22},
      ISSN = {0045-7825,1879-2138},
   MRCLASS = {65N38 (65D17)},
  MRNUMBER = {5025133},
       DOI = {10.1016/j.cma.2026.118773},
       URL = {https://doi.org/10.1016/j.cma.2026.118773},
}

@article{helm,
AUTHOR = {Degli Esposti, Bruno and Falini, Antonella and Kandu{\v c},
              Tadej and Sampoli, Maria Lucia and Sestini, Alessandra},
     TITLE = {Ig{A}-{BEM} for 3{D} {H}elmholtz problems using conforming and
              non-conforming multi-patch discretizations and {B}-spline
              tailored numerical integration},
   JOURNAL = {Comput. Math. Appl.},
  FJOURNAL = {Computers \& Mathematics with Applications. An International
              Journal},
    VOLUME = {147},
      YEAR = {2023},
     PAGES = {164--184},
      ISSN = {0898-1221,1873-7668},
   MRCLASS = {65M60 (74S22)},
  MRNUMBER = {4624891},
       DOI = {10.1016/j.camwa.2023.07.012},
       URL = {https://doi.org/10.1016/j.camwa.2023.07.012},
}

@article {schoenbergwhitney,
    AUTHOR = {Schoenberg, I. J. and Whitney, Anne},
     TITLE = {On {P}\'olya frequence functions. {III}. {T}he positivity of
              translation determinants with an application to the
              interpolation problem by spline curves},
   JOURNAL = {Trans. Amer. Math. Soc.},
  FJOURNAL = {Transactions of the American Mathematical Society},
    VOLUME = {74},
      YEAR = {1953},
     PAGES = {246--259},
      ISSN = {0002-9947,1088-6850},
   MRCLASS = {27.0X},
  MRNUMBER = {53177},
MRREVIEWER = {E.\ Hille},
       DOI = {10.2307/1990881},
       URL = {https://doi.org/10.2307/1990881},
}






\end{document}